\documentclass[onecolumn]{elsarticle}

\journal{Journal of Mathematical Analysis and Applications}

\usepackage[english]{babel}
\usepackage{amsfonts}
\usepackage{amsmath}
\usepackage{amscd}
\usepackage[active]{srcltx} 
\usepackage{latexsym}
\usepackage{amssymb}
\usepackage[latin1]{inputenc}
\usepackage{enumerate}
\usepackage{combelow}

\usepackage{pdfsync}
\usepackage{amsthm}
\usepackage{amsmath}
\usepackage{mathtools}
\usepackage{amscd}
\usepackage{latexsym}
\usepackage{amssymb}
\usepackage{enumitem}

\newtheorem{thm}{Theorem}[section]
\newtheorem{prop}[thm]{Proposition}
\newtheorem{cor}[thm]{Corollary}
\newtheorem*{cor*}{Corollary}
\newtheorem{lema}[thm]{Lemma}
\newtheorem*{lema*}{Lemma}

\numberwithin{equation}{section}
\theoremstyle{definition}
\newtheorem*{Def}{Definition}

\newenvironment{dem}{\vspace{1ex}\noindent{\it Proof.}\hspace{0.5em}}
{\hfill\qed\vspace{1ex}}

\newtheorem{example}{Example}
\newtheorem*{obs}{Remark}

\newtheorem*{obs*}{Remark}
\newtheorem*{thm*}{Theorem}
\newtheorem*{prop*}{Proposition}

\newtheoremstyle{dotless}{}{}{}{}{}{}{ }{}
\theoremstyle{dotless}

\newcommand{\m}[2]{\displaystyle\
    \left[
       \begin{array}{c}
        {#1}\\
        {#2}
       \end{array}
     \right]}

 \newcommand{\ms}[2]{\displaystyle\
	\begin{array}{c}
		{#1}\\
		{#2}
	\end{array}
}

\newcommand{\PI}[2]{\left\langle \,#1 , #2\, \right\rangle}
\newcommand{\pl}{{\mathbin{\!/\mkern-3mu/\!}}}
\newcommand{\ra}{\rightarrow}

\def\minus{\stackrel{-}{\leq}}

\def\lminus{{\;{\vphantom{\leq}}_{-}\leq\,}}

\newcommand{\St}{\mathcal{S}}
\newcommand{\Rt}{\mathcal{R}}
\newcommand{\HH}{\mathcal{H}}
\newcommand{\DD}{\mathcal{D}}
\newcommand{\M}{\mathcal{M}}
\newcommand{\N}{\mathcal{N}}

\newcommand{\Q}{\mathcal{Q}}
\newcommand{\T}{\mathcal{T}}

\newcommand{\G}{\Gamma}
\newcommand{\KK}{\mathcal{K}}
\newcommand{\mc}[1]{\mathcal{#1}}

\newcommand{\PL}{\mc{P} \! \cdot \! L(\HH)^+}
\newcommand{\ol}{\overline}

\begin{document}

\begin{frontmatter}
	
	\title{Semiclosed projections and applications}

	\author[FI,IAM]{Maximiliano Contino\corref{ca}}
	\ead{mcontino@fi.uba.ar}
	
	\author[FI,IAM]{Alejandra Maestripieri}
	\ead{amaestri@fi.uba.ar}
	
		\author[IAM,UNGS]{Stefania Marcantognini}
	\ead{smarcantognini@ungs.edu.ar}

	\cortext[ca]{Corresponding author}
	\address[FI]{%
		Facultad de Ingenier\'{\i}a, Universidad de Buenos Aires\\
		Paseo Col\'on 850 \\
		(1063) Buenos Aires,
		Argentina 
	}
	
	\address[IAM]{%
		Instituto Argentino de Matem\'atica ``Alberto P. Calder\'on'' \\ CONICET\\
		Saavedra 15, Piso 3\\
		(1083) Buenos Aires, 
		Argentina }
	

	\address[UNGS]{%
 Universidad Nacional de General Sarmiento -- Instituto de Ciencias \\ Juan Mar\'ia Gutierrez \\ (1613) Los Polvorines, Pcia. de Buenos Aires, Argentina
}
	
	\begin{abstract}
We characterize the semiclosed projections and apply them to compute the Schur complement of a selfadjoint operator with respect to a closed subspace. These projections occur naturally when dealing with weak complementability.
	\end{abstract}
	
	
\end{frontmatter}

\section{Introduction}
A linear subspace $\mathcal E$ of a Hilbert space $(\HH, \langle\cdot,\cdot\rangle)$ is {\emph{semiclosed}} if there exists an inner product $\langle\cdot,\cdot\rangle^\prime$ such that $(\mathcal E,\langle\cdot,\cdot\rangle^\prime)$ is complete and continuously included in $\HH$. The notion of semiclosed subspace was introduced by Kaufman but subspaces of this type appear in the literature before his seminal paper \cite{Kaufman}. For instance, in the form of {\emph{contractively included subspaces}} in the theory of the de~Branges-Rovnyak spaces \cite{dBR, dBR1}, as the {\emph{operator ranges}} of Fillmore and Williams \cite{Filmore} and, under the name of  {\emph{para-closed}} subspaces, in the work of Foia\cb{s} on the lattice of invariant subspaces \cite{Foias}. Amongst the semiclosed subspaces of $\HH\times \HH$, Kaufman pays much attention to those that are graphs of linear operators in $\HH$, the so-called {\emph{semiclosed}} operators. In fact, the family of all such operators is the main object of analysis in his account of semiclosed subspaces and operators. In the light of Kaufman's study on semiclosed operators we recognize a semiclosed operator $V$ in the factorization $A = BV$ given by the earlier Douglas' Lemma for closed densely defined operators $A, B$ satisfying the operator range inclusion $\Rt(A) \subseteq \Rt(B)$, as well as in the Banach space version of the concept of semiclosed operator in a previous work by Caradus \cite{Caradus}. 

An operator $E$ with domain $\DD(E)$ and range $\Rt(E)$ in $\HH$ is a {\emph{projection}} provided $\Rt(E) \subseteq \DD(E)$ and $E^2x = Ex$  for all $x \in \DD(E)$.
A semiclosed projection $E$ densely defined in a Hilbert space $\HH$ occurs when we are given a closed subspace $\St$ of $\HH$ and a selfadjoint operator $B$ everywhere defined on $\HH$  such that $B$ is $\St$-weakly complementable. Indeed, in this case the Schur complement $B_{/\St}$ of $B$ to $\St$ exists and $B_{/\St} = (I-E)B$ for some (as a matter of fact, any) semiclosed densely defined projection $E$ in an appropriate class.   This fact (cf. \cite{Contino4}) drew our attention to study Hilbert space projections that are densely defined and semiclosed. 

Closed projections were studied by \^Ota in \cite{Ota}, where he showed that a projection is closed if and only if its nullspace and range are both closed subspaces. We generalize the result to  semiclosed projections obtaining the corresponding assertion. Further in-depth 
investigations on closed projections were carried on by Ando \cite{Ando2}. He proved that a closed densely defined projection $E$ with nullspace $\N$  and range $\M$ is distinctively represented as 
$E=(\G^{-1}P_{\M})^*\G^{-1}$ where $\G:=(P_{\M}+P_{\N})^{1/2}$ with $P_\M$ and $P_\N$ the orthogonal projectors onto $\M$ and $\N,$ respectively. Moreover, he established that the well defined operator $\G^{-1}E\G$ is a bounded orthogonal projection. We obtain the analogous result for semiclosed projections with $P_{\M}$ and $P_{\N}$ replaced by $A_1$ and $A_2,$  respectively, where $(A_1, A_2)$ is any pair of positive semidefinite operators such that $\M = \Rt(A_1)$ and $\N = \Rt(A_2)$. However, in this case, it is not possible to obtain a ``distinguished'' representation.

Ando also gave a $2\times 2$ block matrix representation of a closed densely defined projection. The analog for a semiclosed densely defined projection $E$ can be obtained under some extra condition on $\DD(E)$. The very same condition allows us to define the Moore-Penrose pseudoinverse $E^\dagger$. The case when $E$ is closed was considered in \cite{Corach2011}. Therein it was shown that the inverse gives a bijective correspondence between the products of pairs of orthogonal projections and the set of closed densely defined projections. More generally, we prove that the Moore-Penrose pseudoinverse of a semiclosed projection can be related to an operator which is the product of an orthogonal projection times a positive operator. The set of products $PA$, with $P$ an orthogonal projection and $A$ a positive semidefinite operator was studied in \cite{Arias}. 

The semiclosed projections are fundamental when studying weak complementability, a concept introduced in \cite{AntCorSto06} for operators in Hilbert spaces which is a generalization of the notion of complementability, introduced by Ando for matrices \cite{AndoSchur}.
The semiclosed projections arisen in this context, when given a selfadjoint operator $B$ on $\HH$ and a closed subspace $\St$ of $\HH$ such that $B$ is $\St$-weakly complementable, are studied and fully characterized. 
On the other hand, we study the set of  quasi-complementable pairs $(B,\St),$ i.e., the set of $B$-symmetric closed projections onto a prescribed subspace $\St.$ The relation between the notions of weak complementability and quasi-complementability is analyzed to establish whether they are comparable and to what extent. 

Finally, we give a formula of the Schur complement $B_{/ \St}$ of a selfadjoint operator $B$ to $\St$ in terms of semiclosed projections. Also, we characterize $B_{/ \St}$ as the maximum of a set, when a generalization of the minus order is considered, using again semiclosed projections, see \cite{AntCorSto06} and \cite{Arias2}. 

The paper has five sections including this one. Section 2 is a brief expository introduction to
semiclosed subspaces and operators, and serves to set the notation and give some other results
that are needed in the following sections. Section 3 is entirely devoted to the study
of the class of semiclosed densely defined projections. We also deal with the semiclosed densely defined projections having Moore-Penrose inverses, in particular, those with closed nullspaces. In Section 4 we are concerned with $B$-symmetric projections while in Section 5 we study the notions of weak and quasi complementability and give some applications.

\section{Preliminaries}

We assume that all Hilbert spaces are complex and separable. If $\HH$ and $\KK$ are Hilbert spaces, by an \emph{operator} from $\HH$ to $\KK$ we mean a linear function from a subspace of $\HH$ to $\KK.$ The domain, range, nullspace and graph of any given operator $A$ are denoted by $\DD(A),$ $\Rt(A), \ \N(A)$ and $gr(A),$ respectively. 
Given a subset $\T \subseteq \KK,$ the preimage of $\T$ under $A$ is  $A^{-1}(\T):=\{ x \in \HH: \ Ax \in \T \}.$ $L(\HH, \KK)$ stands for the space of the bounded linear operators everywhere defined on $\HH$ to $\KK.$ When $\HH = \KK$ we write, for short, $L(\HH).$ 

The direct sum of two subspaces $\M$ and $\N$ of $\HH$ is represented by $\M \dot{+} \N.$ If, moreover, $\M \perp \N$ their orthogonal sum is denoted by $\M \oplus \N.$ 
The symbol $\Q$ indicates the subset of the oblique projections in  $L(\HH),$ namely, $\Q:=\{Q \in L(\HH): Q^{2}=Q\}$ and  $\mc{P}$ the subset of all the orthogonal projections in $L(\HH),$ $\mc{P}:=\{P \in L(\HH): P^2=P=P^*\};$ for a closed subspace $\M,$ $P_{\M}$ denotes the element in $\mc{P}$ with range $\M.$

Denote by $L(\HH)^s$ the set of selfadjoint operators in $L(\HH),$ $GL(\HH)$ the group of invertible operators in $L(\HH),$ $L(\HH)^+$ the cone of positive semidefinite operators in  $L(\HH)$ and set $GL(\HH)^+:=GL(\HH) \cap L(\HH)^+.$
Given two operators $S, T \in L(\HH),$ the notation  $T \leq S$ signifies that $S-T \in L(\HH)^+.$ 
Given any $T \in L(\HH),$ $\vert T \vert := (T^*T)^{1/2}$ is the modulus of $T$ and $T=U\vert T\vert$ is the polar decomposition of  $T,$ with $U$ the partial isometry such that $\N(U)=\N(T).$ If $T \in L(\HH)^s$ then $T$ has a unique polar decomposition $T=U\vert T\vert$ such that $U=U^*=U^{-1}$ and, moreover, $U$ and $\vert T \vert$ commute.

Given $B \in L(\HH)^s$ and a (non necessarily closed) subspace $\St$ of $\mc{H},$ the $B$-orthogonal complement of $\St$ is $\St^{\perp_{B}}:=\{x \in \mc{H}: \PI{Bx}{y}=0, \mbox{ for every } y\in \St\}=B^{-1}(\St^{\perp})=(B\St)^{\perp}.$

\vspace{0.35cm}
The next result, due to Fillmore and Williams, characterizes the sum and the intersection of operator ranges as operator ranges.

\begin{thm}[{\cite[Theorem 2.2, Corollary 2]{Filmore}}]  \label{FW} Let $A, B \in L(\HH).$ Then
\begin{enumerate}
	\item $\Rt(A)+\Rt(B)=\Rt((AA^*+BB^*)^{1/2}).$
	\item There exist $X, Y \in L(\HH)$ such that $\Rt(A) \cap \Rt(B) = \Rt((AXA^*)^{1/2})=\Rt((BYB^*)^{1/2}).$
\end{enumerate}
\end{thm}

Given $B \in L(\HH)^s$ and $\St$ a closed subspace of $\HH,$ we say that $\St$ is $B$-\emph{positive} if 
$\PI{Bs}{s} > 0 \mbox{ for every } s \in \St, \ s\not =0.$ $B$-\emph{nonnegative}, $B$-\emph{neutral}, $B$-\emph{negative} and $B$-\emph{nonpositive} subspaces are defined analogously. If $\St$ and $\T$ are two closed subspaces of $\HH,$ the notation $\St \ \oplus_{B} \  \T$ is used to indicate the direct sum of $\St$ and $\T$ when, in addition, $\PI{Bs}{t}=0 \mbox{ for every } s \in \St \mbox{ and } t \in \T.$ 

The following is a consequence of the spectral theorem for Hilbert space selfadjoint operators. 

\begin{lema} \label{lemaBdecom} Let $B \in L(\HH)^s$ and $\St$ be a closed subspace of $\HH.$  Then the Grammian of $B,$ $G_{B,\St}:=P_{\St}B|_\St$ can be represented as
	\begin{equation} \label{Bdecomp}
		G_{B,\St}=G_1-G_2,\\
	\end{equation}
	where $G_1, G_2 \in L(\St)^+$ and $\Rt(G_1) \perp \Rt(G_2).$ 
	Also, if $\St_+:=\ol{\Rt(G_1)}$ and $\St_-:=\ol{\Rt(G_2)}\oplus  \N(G_{B,\St}),$ then $\St_+ \perp \St_-$ and $\St$ can be represented as
	\begin{equation} \label{Sdecomp}
		\St = \St_{+} \ \oplus_{B} \ \St_{-},
	\end{equation}
	where $\St_+$ is $B$-positive and $\St_-$ is $B$-nonpositive.
\end{lema}

Let $A, B \in L(\HH).$ Then, by Douglas' Lemma \cite{Douglas}, $\Rt(B) \subseteq \Rt(A)$ if and only the equation $AX=B$ has a solution in $L(\HH).$ In this case, there exists a unique $D$ such that $AD=B$ and $\Rt(D) \subseteq \ol{\Rt(A^*)}.$ The operator $D$ is called the \emph{reduced solution} of the equation $AX=B.$

The next lemma characterizes the positive operators in terms of its matrix decomposition, see \cite{Shorted2}.

\begin{lema} \label{LemmaPositive} Let $\St \subseteq \HH$ be a closed subspace and $B \in L(\HH)^s$  with matrix decomposition $$B=\begin{bmatrix}
	a  & b \\ 
	b^* & c \\
	\end{bmatrix} \ms{\St}{\St^{\perp}}.$$ Then $B \in L(\HH)^+$ if and only if 
	$$a \geq 0,  \, \,  \Rt(b) \subseteq \Rt(a^{1/2}) \mbox{ and } c=f^*f+t,$$ for some $t \geq 0,$ where $f$ is the reduced solution of the equation $b=a^{1/2}x.$
\end{lema}

\subsection*{\textbf{Semiclosed subspaces and  operators}}
The notions of \emph{semiclosed subspace} and \emph{semiclosed operator} were formally introduced by Kaufman \cite{Kaufman}, though these notions were considered by other authors before, as we pointed out in the Introduction.

\begin{Def}
A subspace $\St$ of $\HH$ is \emph{semiclosed} if $\St$ is a (not necessarily closed) subspace for which there exists an inner product $\PI{\cdot}{\cdot}^{'}$ such that $(\St,\PI{\cdot}{\cdot}^{'})$ is a Hilbert space which is \emph{continuosly included} in $\HH,$ i.e., there exists $b > 0$ such that $\PI{x}{x} \leq b \PI{x}{x}'$ for every $x \in \St.$ 
\end{Def}
As only an infinite dimensional subspace can be semiclosed but not closed, only infinite dimensional complex Hilbert spaces are considered. 

Operator ranges are semiclosed subspaces: in fact, if $T \in L(\HH)$ define
$$\Vert u \Vert_T := \Vert T^{\dagger} u \Vert \mbox{ for } u \in \Rt(T),$$ where $T^{\dagger}$ denotes the (possibly unbounded) Moore-Penrose inverse of $T,$ see \cite{Nashed}.
Then $(\Rt(T), \Vert \cdot \Vert_T)$ is a Hilbert space and
\begin{equation} \label{eqT}
\Vert u \Vert = \Vert TT^{\dagger} u \Vert \leq \Vert T \Vert \Vert T^{\dagger} u \Vert =\Vert T \Vert \Vert  u \Vert_T \mbox{ for } u \in \Rt(T).
\end{equation}
See \cite{Ando, Bolotnikov}.

The space $\Rt(T)$ equipped with the Hilbert space structure $\Vert \cdot \Vert_T$ is denoted by $$\mc{M}(T) := (\Rt(T), \Vert \cdot \Vert_T).$$ The Hilbert spaces $\mc{M}(T)$ play a significant role in many areas, in particular in the de Branges complementation theory \cite{Ando}. 

The semiclosed subspaces are all of them operator ranges: Fillmore and Williams proved that  $\St$ is a semiclosed subspace of $\HH$ if and only if $\St$ is the range of a closed operator $T$ on $\HH.$ Moreover, the operator $T$ can be chosen to be bounded and positive (semidefinite), see  \cite[Theorem 1.1]{Filmore}.
Furthermore, if $T$ is a contraction, i.e. $\Vert T \Vert \leq 1,$ then $\St':=\mc{M}((I-TT^*)^{1/2})$ is its \emph{de Branges complement} and  $\St + \St'=\HH$ \cite[Corollary 3.8]{Ando}, where the last sum need not be direct \cite[Proposition 3.4]{Bolotnikov}.

Given two operators $T_1, T_2 \in L(\HH),$ by Theorem \ref{FW}, the subspace $\Rt(T_1)+\Rt(T_2)$ is the range of $T:=(T_1T_1^{*}+T_2T_2^{*})^{1/2}.$ This shows that the sum of semiclosed subspaces is again semiclosed.
The following interesting result by Ando compares the norm $\Vert \cdot \Vert_T$  with the norms $\Vert \cdot \Vert_{T_1}$ and $\Vert \cdot \Vert_{T_2}.$

\begin{thm}[{\cite[Corollary 3.8]{Ando}}] \label{thmAndo} For $T_1, T_2 \in L(\HH),$ let $T:=(T_1T_1^*+T_2T_2^*)^{1/2}.$ Then $\Vert u_1 + u_2 \Vert_T^2 \leq \Vert u_1 \Vert_{T_1}^2 + \Vert u_2 \Vert_{T_2}^2,$ for $u_1 \in \Rt(T_1)$ and $u_2 \in \Rt(T_2),$ and for any $u \in \Rt(T),$ there are unique  $u_1 \in \Rt(T_1)$ and $u_2 \in \Rt(T_2)$ such that $u=u_1+u_2$ and $$ \Vert u_1 + u_2 \Vert_T^2 = \Vert u_1 \Vert_{T_1}^2 + \Vert u_2 \Vert_{T_2}^2.$$
\end{thm}
Applying again Theorem \ref{FW} it follows that the family of semiclosed subspaces is closed under intersection, see also \cite[Proposition 4, Proposition 6]{Balaji}. The set of semiclosed subspaces is the lattice of domains of closed operators in $\HH$ \cite{Filmore}. Also, if $\St$ is a semiclosed subspace of $\HH$ then all inner products $\PI{\cdot}{\cdot}^{'}$ such that $(\St,\PI{\cdot}{\cdot}^{'})$ is a Hilbert space which is continuosly included in $\HH$ generate the same topology on $\St.$ See \cite{Macnearney1, Macnearney2} and \cite[Theorem 11]{Balaji}.

\bigskip 
\begin{Def} [\cite{Kaufman,Koliha}] An operator $C: \DD(C) \subseteq \HH \ra \KK$ is a \emph{semiclosed operator} if $gr(C)$ is a semiclosed subspace of $\HH \times \KK.$
\end{Def}

Denote by $SC(\HH, \KK)$ the set of all semiclosed operators with domain in $\HH$ to $\KK$ and set $SC(\HH):=SC(\HH,\HH).$ The following is a characterization of $SC(\HH),$ see \cite[Theorem 1]{Kaufman}. 

\begin{thm} \label{thmKaufman} Given an operator $C : \DD(C)\subseteq \HH \ra \HH,$ the following are equivalent:
	\begin{itemize}
		\item [i)] $C \in SC(\HH);$
		\item [ii)] $\DD(C)$ is a semiclosed subspace of $\HH$ and $C \in L(\DD(C),\HH);$
		\item [iii)] there exist $A \in L(\HH)$ and  $D \in L(\HH)^+$ such that $C=AD^{\dagger}|_{\Rt(D)}$ and $\N(D) \subseteq \N(A).$
	\end{itemize}
\end{thm}

\begin{cor} \label{corKaufman} Let  $C: \DD(C) \subseteq \HH \ra \HH$ be a given operator. Then $C \in SC(\HH)$ if and only if there exists $D \in L(\HH)^+$ such that $\DD(C)=\Rt(D)$ and $CD \in L(\HH).$
\end{cor}

\begin{cor} \label{corKaufman2} Let  $C \in SC(\HH).$ Then $CA  \in L(\KK,\HH)$ for any $A \in L(\KK,\HH)$ such that $\Rt(A) \subseteq \DD(C).$
\end{cor}
\begin{dem} By Corollary \ref{corKaufman},  there exists an operator $D \in L(\HH)^+$ such that $\DD(C)=\Rt(D)$ and $CD \in L(\HH).$ Then, if $\Rt(A) \subseteq \DD(C)=\Rt(D),$ by Douglas' Lemma, $A=DX_0$ for some $X_0 \in L(\KK, \HH).$ Therefore, $CA=CDX_0 \in L(\KK,\HH).$ 
\end{dem}

The set $SC(\HH)$ is closed under addition, multiplication, inversion and restriction to semiclosed subspaces of $\HH$  \cite{Kaufman}. Also, if $T_1,T_2 \in SC(\HH,\KK)$ are such that $T_1$ and $T_2$ coincide on $\DD(T_1) \cap \DD(T_2),$ then the operator $T : \DD(T_1) + \DD(T_2) \ra \KK$ coinciding with $T_1$ on $\DD(T_1)$ and with $T_2$ on $\DD(T_2)$ is a semiclosed operator \cite{Djikic}.

\begin{obs} In \cite[Theorem 2]{Douglas}, Douglas proved that given $A, B$ densely defined closed operators on $\HH$ such that $\Rt(A) \subseteq \Rt(B),$  there exist an operator $V$ on $\HH$ with $\DD(V)=\DD(A)$ and a number $M \geq 0$ such that 
\begin{equation} \label{DouglasEq}
A=BV \mbox{ and } \Vert Vx \Vert^2 \leq M (\Vert x \Vert^2 + \Vert Ax \Vert^2), \ \mbox{ for every } x \in \DD(V).
\end{equation}

The operator $V$ is semiclosed: in fact, define $$\PI{(x,Vx)}{(y,Vy)}':=\PI{x}{y}+\PI{Ax}{Ay} \mbox{ for } x, y \in \DD(V).$$ Then $(gr(V), \PI{\cdot}{\cdot}')$ is a Hilbert space because $A$ is closed. On the other hand, since \eqref{DouglasEq} holds, 
$$ \PI{(x,Vx)}{(x,Vx)} \leq (M+1) \PI{(x,Vx)}{(x,Vx)}' \mbox{ for every } x \in \DD(V).$$ Hence $(gr(V), \PI{\cdot}{\cdot}')$ is continuously included in $\HH \times \HH.$

\end{obs}
\section{Semiclosed projections}

A linear operator $E$ acting in $\HH$ is a \emph{projection} if
$$\Rt(E) \subseteq \DD(E) \mbox{ and } E^2x=Ex \mbox{ for every } x \in \DD(E).$$

\begin{thm}[{\cite[Lemma 3.5]{Ota}}] \label{ProjectionsOTA}
\begin{enumerate}
\item If $E$ is a projection in $\HH$ then $$\Rt(E) \dotplus \N(E)= \DD(E).$$ Conversely, given two subspaces ${\N, \M}$ of $\HH$ such that $\N \cap\M = \{0\},$ there exists a projection $E$ with $\Rt(E)=\M$ and $\N(E)=\N.$
\item $E$ is a closed projection if and only if $\Rt(E)$ and $\N(E)$ are closed subspaces of $\HH.$ 
\end{enumerate} 
\end{thm}

Write $E=P_{\M \pl \N}$ to denote the projection with $\Rt(E)=\M$ and $\N(E)=\N.$
If $E$ is a densely defined projection in $\HH,$ then $E^*$ is a (non necessarily densely defined) closed projection (see \cite[Proposition 3.4]{Ota}) with $\N(E^*)=\Rt(E)^{\perp}$ and $\Rt(E^*)=\N(E)^{\perp}.$ The last equality follows from the former and the fact that $I-E$ is a projection with domain $\DD(E)$ so that $\Rt(E^*)=\N(I-E^*)=\Rt(I-E)^{\perp}=\N(E)^{\perp}.$
Then, by Theorem \ref{ProjectionsOTA}, $$\DD(E^*)=\N(E)^{\perp} \dotplus \Rt(E)^{\perp}.$$
From this, it is immediate that if $E$ is a densely defined closed projection, then $E^*$ is densely defined.

Every densely defined closed projection $E=P_{\M \pl \N}$ admits a matrix representation according to the decomposition $\HH=\M \oplus \M^{\perp},$
\begin{equation} \label{matrixE}
	E=\begin{bmatrix}
		1&x\\
		0& 0
	\end{bmatrix}
\end{equation} such that $x : \DD(x) \subseteq \M^{\perp} \ra \M$ is a densely defined closed operator, where $\DD(x)=P_{\M^{\perp}}(\N).$ 

Conversely, consider the decomposition $\HH=\M \oplus \M^{\perp}.$ If
\begin{equation} \label{matrixE2}
E=\begin{bmatrix}
1&x\\
0& 0
\end{bmatrix},
\end{equation} where $x : \DD(x) \subseteq \M^{\perp} \ra \M$ is a densely defined operator and $\DD(E)=\M \oplus \DD(x).$ Then $E$ is a densely defined projection with $\Rt(E)=\M.$ Furthermore, $E$ is closed  if and only if $x$ is closed. In this case, 
\begin{equation} \label{Esclosed}
E^*=\begin{bmatrix}
	1&0\\
	x^*& 0
\end{bmatrix},
\end{equation}
where $x^*: \DD(x^*) \subseteq \M \ra \M^{\perp}$ is densely defined. See \cite[Proposition 1.7]{Cuasi} for the proof of these assertions.

\bigskip 
More generally,
\begin{prop}\label{semiclosedProjections} Let $E: \DD(E) \subseteq \HH \ra \HH$ be a projection. Then $E$ is semiclosed if and only if $\Rt(E)$ and $\N(E)$ are semiclosed subspaces.
\end{prop}

\begin{dem} If $\Rt(E)$ and $\N(E)$ are semiclosed subspaces of $\HH,$  there exist $A_1, A_2 \in L(\HH)^+$ such that $\Rt(E)=\Rt(A_1)$ and $\N(E)=\Rt(A_2).$ Then, by Theorem \ref{FW}, $\DD(E)=\Rt(E)\dotplus \N(E)$ is semiclosed. Let us see that $E \in L(\DD(E),\HH).$   
Consider $\G=(A_1^2+A_2^2)^{1/2}.$ Then, by Theorem \ref{FW}, $\Rt(\G)=\Rt(A_1)\dotplus\Rt(A_2)=\DD(E).$ 
Let $u \in \DD(E).$ Then, by Theorem \ref{thmAndo}, there exist uniquely $m \in \Rt(E)$ and $n \in \N(E)$ such that $u=m+n$ and $ \Vert u \Vert_\G^2 = \Vert m \Vert_{A_1}^2 + \Vert n \Vert_{A_2}^2.$ Then, using \eqref{eqT}, $$\Vert Eu\Vert^2 = \Vert m \Vert^2 \leq \Vert A_1 \Vert^2 \Vert m \Vert_{A_1}^2 \leq \Vert A_1 \Vert^2 (\Vert m \Vert_{A_1}^2 +\Vert n \Vert_{A_2}^2) =\Vert A_1 \Vert^2  \Vert u \Vert_\G^2.$$ Then $E \in L(\DD(E),\HH)$ and, by Theorem \ref{thmKaufman}, $E \in SC(\HH).$

Conversely, suppose that $E$ is a semiclosed projection. Then, by Theorem \ref{thmKaufman}, $\DD(E)$ is a semiclosed subspace of $\HH$ and, by \cite[Theorem 2]{Kaufman}, $\Rt(E)=E(\DD(E))$ is also a semiclosed subspace of $\HH.$ Since the set $SC(\HH)$ is closed under addition, $E \in SC(\HH)$ if and only if $I-E \in SC(\HH).$ 
Hence, $\N(E)=\Rt(I-E)$ is a semiclosed subspace of $\HH,$ where we used again \cite[Theorem 2]{Kaufman}. 
\end{dem}

\begin{Def} Let $$\Q_{SC}:=\{E: E \mbox{ is a densely defined semiclosed projection on } \HH\}.$$
\end{Def}
\bigskip
In \cite[Theorem 2.2]{Ando2}, Ando proved that if $E=P_{\M \pl \N}$ is a closed projection and $\G:=(P_{\M}+P_{\N})^{1/2}$ then $\DD(E)=\Rt(\G)$ and $E$ admits the following representation:
$$E=(\G^{-1}P_{\M})^*\G^{-1}.$$ 
Moreover, the well defined operator $P:=\G^{-1}E\G$ is an orthogonal projection, see \cite[Theorem 2.3]{Ando2}.
Analogous results can be obtained for densely defined semiclosed projections. 

Let $E=P_{\M \pl \N} \in \Q_{SC}.$  Since, by Proposition \ref{semiclosedProjections}, $\M$ and $\N$ are semiclosed subspaces, there exist $A_1, A_2 \in L(\HH)^+$ such that $\M=\Rt(A_1)$ and $\N=\Rt(A_2).$ Define the operator $\G=\G(A_1,A_2)$ as
\begin{equation} \label{Gamma}
\G:=(A_1^2+A_2^2)^{1/2}. 
\end{equation}
Then  $\G \in L(\HH)^+$ and, by Theorem \ref{FW}, $\Rt(\G)=\Rt(A_1)\dotplus \Rt(A_2)=\DD(E)$ is dense, so that $\G$ is  injective.

\begin{prop}[{cf. \cite[Theorem 2.2]{Ando2}}] \label{ddprojection2} Let $E\in \Q_{SC}$ with $\Rt(E)=\Rt(A_1),$ $\N(E)=\Rt(A_2),$ $A_1, A_2 \in L(\HH)^+$ and $\G$ as in \eqref{Gamma}. Then $E$ admits the representation 
$$E=(\G^{-1}A_1^2)^*\G^{-1}.$$
\end{prop}
\begin{dem} Since $\Rt(A_1^2) \subseteq \Rt(A_1) \subseteq \Rt(\G)$ and $\Gamma$ is injective, by Douglas' Lemma, there exists a unique $D \in L(\HH)$ such that $A_1^2= \G D^*=D \G,$ then $D^*=\G^{-1}A_1^2$ and $D=(\G^{-1}A_1^2)^*.$
Write $\tilde{E}=D \G^{-1}.$ Since $\DD(E)=\Rt(\G)=\DD(\tilde{E}),$ for the proof of the assertion it suffices to show that $E\G x=\tilde{E}\G x=Dx \mbox{ for every } x \in \HH.$ Since $E\G$ is bounded  by Corollary \ref{corKaufman2}, and $\Rt(\G)$ is dense in $\HH,$ the equality is guaranteed if the operators $E\G$ and $D$ coincide on this dense subspace. It is
$E\G^2 x= E(A_1^2+A_2^2)x=A_1^2x=D\G x.$ Therefore $E\G=D$ or $E=\tilde{E}=D\G^{-1}$ on $\Rt(\G).$ 
\end{dem}

\begin{cor}[{cf. \cite[Theorem 2.3]{Ando2}}] \label{ddprojection} Let $E$ be a densely defined operator in $\HH.$ Then $E \in \Q_{SC}$ if and only if there exists $\G \in L(\HH)^+$ injective with $\Rt(\G)=\DD(E)$ such that $$\G^{-1}E\G \in \mc{P}.$$
\end{cor}
\begin{dem} 
Suppose that $E\in \Q_{SC}$ with $\Rt(E)=\Rt(A_1),$ $\N(E)=\Rt(A_2),$ $A_1, A_2 \in L(\HH)^+.$ Let $D:=(\G^{-1}A_1^2)^*,$ with $\G$ as in \eqref{Gamma}.  Then $\G$ is injective and, by Proposition \ref{ddprojection2}, $E\G=D.$ Let $P_\G:=\G^{-1}E\G=\G^{-1}D,$ then $P_\G$ is a bounded projection. In fact, since  $\Rt(D) \subseteq \Rt(E) \subseteq \DD(E) = \Rt(\G),$ by Douglas' Lemma, the  only solution $X_0$ of the equation $D=\G X$ is given by $X_0=\G^{-1}D \in L(\HH).$ Also, since $\Rt(A_1) \subseteq \Rt(\G),$ by Douglas's Lemma again, there exists a unique $D' \in L(\HH)$ such that $A_1= \G D'^*=D' \G,$ then $D'^*=\G^{-1}A_1$ and $D'=(\G^{-1}A_1)^*.$ 
From $\G (D'^*\G D'^*)=A_1^2=\G D^*$ and the fact that $\G$ is injective, it follows that $D=D' \G D'.$ Then $P_\G=\G^{-1} D=\G^{-1} D' \G D'=\G^{-1}A_1D'=(D')^*D' \in L(\HH)^s.$ Also $P_\G^2=\G^{-1}D\G^{-1}D=\G^{-1}E\G\G^{-1}D=\G^{-1}ED=\G^{-1}D=P_\G.$

Conversely, suppose that $E$ is a densely defined operator in $\HH$ such that there exists $\G \in L(\HH)^+$ with $\Rt(\G)=\DD(E)$ and $\G^{-1}E\G:=P_\G \in \mc{P}.$ Then $E\G=\G P_\G$ and, by Corollary \ref{corKaufman}, $E \in  SC(\HH).$ Since $\G$ is injective, $E=\G P_\G\G^{-1},$ $\Rt(E) \subseteq \DD(E)$ and $E^2=E.$ 
\end{dem}

Since $E \G=\G P_\G,$ it follows that $P_\G=P_{\G^{-1}(\Rt(E))}.$ Since $\Rt(E)=\G (\Rt(P_\G)),$ $\ \N(E)= \G(\N(P_\G))$ and $P_\G$ is orthogonal, it holds that $\G^{-1}(\Rt(E)) \perp \G^{-1}(\N(E)).$

If $\G, \G' \in L(\HH)^+$ are as in Corollary \ref{ddprojection}, then $\Rt(\G)=\DD(E)=\Rt(\G')$ and, by Douglas' Lemma, there exists $G \in GL(\HH)$ such that $\G'=\G G=G^*\G.$ 
Moreover, if $P_\G=P_{\G^{-1}(\Rt(E))}$ and $P_{\G'}=P_{\G'^{-1}(\Rt(E))}$ then $$P_{\G'}=G^{-1}P_\G G.$$ In fact, the projection $G^{-1}P_\G G$ is bounded, $\Rt(G^{-1}P_\G G)=G^{-1}(\Rt(P_\G))=(\G G)^{-1}(\Rt(E))=\G'^{-1}(\Rt(E))\linebreak =\Rt(P_{\G'})$ and $\N(G^{-1}P_\G G)=G^{-1}(\N(P_\G))=(\G G)^{-1}(\N(E))=\G'^{-1}(\N(E))=\N(P_{\G'}).$

\subsection{On the Moore-Penrose inverse of semiclosed projections with closed nullspace}

In order to define the Moore-Penrose inverse of a densely defined projection $E$ in a satisfactory fashion we need an extra condition on its domain. This condition guarantees the existence of an orthogonal complement of $\N(E)$ relative to $\DD(E).$

\begin{lema} \label{lemaSP} Let $E=P_{\M \pl \N}$ be a densely defined projection. Then the following statements are equivalent:
	\begin{enumerate}
		\item [i)] $\M \subseteq \N^{\perp} \oplus \N;$
		\item [ii)] $\DD(E)=P_{\N^{\perp}}(\M)\oplus \N;$
		\item [iii)] $\DD(E)=\DD(E) \cap \N^{\perp} \oplus \N.$
	\end{enumerate}

In this case, $\M \cap \ol{\N}=\{0\}$ and therefore $E$ admits the extension $P_{\M \pl \ol{\N}}.$
\end{lema}
\begin{dem} $i) \Rightarrow ii):$ Let $x \in \M.$ Then $x=P_{\N^{\perp}}x+n,$ for some $n \in \N.$ Therefore $x \in P_{\N^{\perp}}(\M)\oplus \N,$ $\M \subseteq  P_{\N^{\perp}}(\M)\oplus \N$ and $\DD(E)=\M \dotplus \N \subseteq P_{\N^{\perp}}(\M)\oplus \N.$ To see the other inclusion, if $y \in P_{\N^{\perp}}(\M)$ then there exists $m \in \M$ such that $y=P_{\N^{\perp}}m.$ Since $m \in \M,$ there exists $t\in \N^{\perp}$ and  $n \in \N$ such that $m=t+n.$ Then
$y=m-(I-P_{\N^{\perp}})m=m-n \in \M \dotplus \N.$ Then $P_{\N^{\perp}}(\M) \subseteq \DD(E)$ and $P_{\N^{\perp}}(\M) \oplus \N \subseteq \DD(E).$ 

$ii) \Rightarrow iii):$ Clearly, $P_{\N^{\perp}}(\M) \subseteq \DD(E) \cap \N^{\perp}.$ On the other hand, let $x \in \DD(E) \cap \N^{\perp}.$ Then $x=m+n,$ for some $m\in \M$ and $n \in \N$ and $x=P_{\N^{\perp}}x=P_{\N^{\perp}}m \in P_{\N^{\perp}}(\M).$  Therefore $P_{\N^{\perp}}(\M)=\DD(E) \cap \N^{\perp}$ and
$\DD(E)=\DD(E) \cap \N^{\perp} \oplus \N.$

$iii) \Rightarrow i):$ It follows from the fact that $\M \subseteq \DD(E)=\DD(E) \cap \N^{\perp} \oplus \N \subseteq \N^{\perp} \oplus \N.$

In this case, from  $\M \subseteq \N \oplus \N^{\perp},$ we have that $P_{\ol{\N}}(\M)\subseteq \N.$ Therefore $\M \cap \ol{\N}=\{0\}.$ In fact, if $x \in \M \cap \ol{\N}.$ Then $x=P_{\ol{\N}}x \in \N \cap \M=\{0\}.$  
\end{dem}

Given $E$ a densely defined projection, condition $i)$ on $\Rt(E)$ and $\N(E)$ need not hold: the example on page 278 of \cite{Filmore} shows that there exist subspaces $\M,\N$ of $\HH$ such that $\ol{\M\dotplus \N}=\HH,$ but $\M \cap \ol{\N} \not = \{0\}.$

\bigskip
We begin by applying Lemma \ref{lemaSP} to obtain a matrix decomposition of a given projection.
Every closed projection $P_{\M  \pl\N}$ admits a matrix representation  according to the orthogonal decomposition $\HH=\N^{\perp} \oplus \N$ on $\DD(E)=P_{\N^{\perp}}(\M)\oplus \N,$ see \cite[Theorem 2.6]{Ando2}. We generalize this result for a densely defined projection $E=P_{\M \pl \N}$  satisfying the conditions of Lemma \ref{lemaSP}.

\begin{prop}[{cf. \cite[Theorem 2.6]{Ando2}}] \label{proprepres} Let $E=P_{\M \pl \N}$ be a densely defined  projection such that $\M \subseteq \N \oplus \N^{\perp}.$ According to the orthogonal decomposition $\HH=\N^{\perp} \oplus \ol{\N},$ $E$ admits the matrix representation 
$$E=\begin{bmatrix}
	I & 0 \\ 
	P_{\ol{\N}}(P_{\N^{\perp}}P_{\ol{\M}}|_\M)^{-1} & 0 \\
\end{bmatrix} 
 \mbox{ on }  P_{\N^{\perp}}(\M) \oplus \N.$$
\end{prop}

\begin{dem} By Lemma \ref{lemaSP},  the operator $P_{\N^{\perp}}P_{\ol{\M}}|_\M$ is injective with dense range in $\N^{\perp}.$ In fact, if $P_{\N^{\perp}}P_{\ol{\M}}x=0$ for $x \in \M,$ then $x \in \M \cap \ol{\N}=\{0\} .$ Also,  $\N^{\perp}=P_{\N^{\perp}}(\ol{\M\dotplus \N}) \subseteq \ol{P_{\N^{\perp}}(\M)}\subseteq \N^{\perp}.$ So that $\ol{\Rt(P_{\N^{\perp}}P_{\ol{\M}}|_\M)}=\ol{P_{\N^{\perp}}(\M)}=\N^{\perp}$ and the operator $P_{\ol{\N}}(P_{\N^{\perp}}P_{\ol{\M}}|_\M)^{-1}$ is a linear operator from $P_{\N^{\perp}}(\M) \subseteq \N^{\perp}$ to $\N.$
	
By Lemma \ref{lemaSP} again, $\DD(E)=P_{\N^{\perp}}(\M)\oplus \N$ and $P_{\N^{\perp}} (\M)=\DD(E) \cap \N^{\perp}.$ Clearly, $EP_{\ol{\N}}|_{\N}=0$ and $P_{\N^{\perp}}EP_{\N^{\perp}}|_{\DD(E) \cap \N^{\perp}}=I_{\N^{\perp}}|_{\DD(E) \cap \N^{\perp}}.$ Finally, $P_{\ol{\N}}EP_{\N^{\perp}}|_{\DD(E) \cap \N^{\perp}}=	P_{\ol{\N}}(P_{\N^{\perp}}P_{\ol{\M}}|_\M)^{-1}$ on \; $P_{\N^{\perp}}(\M).$
\end{dem}

Let $E=P_{\M \pl \N}$ be a densely defined  projection such that $ \M \subseteq \N^{\perp} \oplus \N.$ By Lemma \ref{lemaSP}, $\DD(E) \cap \N^{\perp}=P_{\N^{\perp}}(\M).$ In this case, the Moore-Penrose inverse $E^{\dagger}$ of $E$ is well defined (see \cite{Nashed}):
$E^{\dagger} : \M \oplus \M^{\perp} \ra P_{\N^{\perp}}(\M) \subseteq P_{\N^{\perp}}(\M)\oplus \ol{\N},$ 

$$E^{\dagger} = \left\{
\begin{array}{@{}l@{\thinspace}l}
	0  & \ \ \text{if } x \in \M^{\perp}\\
	(E|_{P_{\N^{\perp}}(\M)})^{-1}x  & \ \ \text{if } x \in \M.\\
\end{array}
\right.
$$

The operators $EE^{\dagger}$ and $E^{\dagger}E$ are well defined and they are densely defined projections. In fact, $EE^{\dagger}=P_{\ol{\M}} \mbox{ on } \DD(E^{\dagger})=\M \oplus \M^{\perp}$ and $E^{\dagger}E=P_{\N^{\perp}} \mbox{ on } \DD(E)=P_{\N^{\perp}}(\M)\oplus \N.$

\begin{prop} \label{corextension} Let $E=P_{\M \pl \N}$ be a densely defined projection such that $\M \subseteq \N^{\perp} \oplus \N.$ Then 
	$$E^{\dagger}=({P_{\M \pl \ol{\N}}})^{\dagger}.$$
\end{prop}

\begin{dem} By Lemma \ref{lemaSP}, $\tilde{E}:=P_{\M \pl \ol{\N}}$ is a densely defined projection such that $E \subseteq \tilde{E}.$ Since $E|_{P_{\N^{\perp}}(\M)}=\tilde{E}|_{P_{\N^{\perp}}(\M)},$ it follows that  $E^{\dagger}={\tilde E}^{\dagger}.$
\end{dem}

In view of Proposition \ref{corextension}, we focus our attention on the Moore-Penrose inverse of densely defined projections with closed nullspace.
\begin{prop} \label{MPSC} Let $E=P_{\M \pl \N}$ be a densely defined  projection with closed nullspace. Then $$E^{\dagger}=P_{\N^{\perp}} P_{\ol{\M}} \mbox{ on } \DD(E^{\dagger})=\M \oplus \M^{\perp}.$$
\end{prop}
\begin{dem} If $x \in \M^{\perp}$ then $E^{\dagger}x=0.$ 
On the other hand, if $m \in \M,$ let $y :=(E|_{P_{\N^{\perp}}(\M)})^{-1}m,$ then $y \in  P_{\N^{\perp}}(\M)$ and $Ey=m;$ also, $EP_{\N^{\perp}}m=Em=m.$ Therefore
$$E^{\dagger}m=(E|_{P_{\N^{\perp}}(\M)})^{-1}m=P_{\N^{\perp}}m=P_{\N^{\perp}}P_{\ol{\M}}m.$$
\end{dem}

If $E=P_{\M \pl \N}$ is a densely defined closed projection, then $E^{\dagger}=P_{\N^{\perp}}P_{\M} \in \mc{P} \cdot \mc{P}.$ Moreover, the map $E \mapsto E^{\dagger}$ from the set of densely defined closed projections onto $\mc{P} \cdot \mc{P}$ is a bijection, see \cite{Corach2011}. 

To study the semiclosed case, consider the set $$\mc{P} \! \cdot \! L(\HH)^+:=\{ T \in L(\HH): T=PA \mbox{ with } P \in \mc{P} \mbox{ and } A \in L(\HH)^+ \}.$$ 
This set was studied in \cite{Arias}, where it was showed that any $T \in \PL$ can be factored as $T=P_TA,$ where $P_T:=P_{\ol{\Rt(T)}}$  and $A \in L(\HH)^+$ is such that $\N(T)=\N(A),$ though this factorization may not be unique. We say that $A \in L(\HH)^+$ is \emph{optimal} for $T$ if $T=P_TA$ and $\N(T)=\N(A).$
A description of the set of optimal operators for $T$ can be found in \cite[Remark 4.2 and Proposition 4.4]{Arias}.

\begin{prop} \label{propMP}
Let $E=P_{\M \pl \N}\in \Q_{SC}.$ If $\N(E)$ is closed then there exists $\G \in L(\HH)^+$ such that $\Rt(\G)=\DD(E^{\dagger})$ and $E^{\dagger} \G \in \PL.$
\end{prop}

\begin{dem} By Proposition \ref{MPSC}, $E^{\dagger}=P_{\N^{\perp}} P_{\ol{\M}}$ $\mbox{ on } \DD(E^{\dagger}).$
If $A \in L(\HH)^+$ is such that $\Rt(A)=\M,$  take $\G:=\begin{bmatrix}
A & 0 \\ 
0 & I \\
\end{bmatrix} \ms{\ol{\M}}{\M^{\perp}}.$ Then $\G \in L(\HH)^+$ and $\Rt(\G)=\Rt(A)\oplus \M^{\perp}=\DD(E^{\dagger}).$
Hence $E^{\dagger}\G=P_{\N^{\perp}} P_{\ol{\M}}\G=P_{\N^{\perp}} A \in \PL.$
\end{dem}

This generalizes the fact that if $E$ is a densely defined closed projection then $E^{\dagger} \in \mc{P} \cdot \mc{P},$ since in this case the operator $\G$ can be chosen to be the identity. From Proposition \ref{propMP} it follows that every $E\in \Q_{SC}$ with closed nullspace has an associated set in $\PL,$ namely, $\{T=P_{\N(E)^{\perp}}A: \ A \in L(\HH)^+, \ \Rt(A)=\Rt(E)\}.$

On the other hand, every $T \in \PL$ has an associated set of semiclosed projections. 
\begin{prop} \label{propMP2}Let $T \in \PL.$ If $A \in L(\HH)^+$ is optimal for $T$ then $$\ol{\Rt(A) \dotplus \Rt(T)^{\perp}}=\HH.$$
\end{prop}
\begin{dem}
Consider $T \in \PL$ and write $T=P_{T}A$ with $A \in L(\HH)^+$  optimal. Observe that $\Rt(A) \cap \Rt(T)^{\perp}=\{0\}.$ In fact, if $x \in \Rt(A) \cap \Rt(T)^{\perp}$ then $x=A y$ for some $y \in \HH$ and $0=P_Tx=P_TAy=Ty.$ Then $y \in \N(T)=\N(A)$ and $x=Ay=0.$ Also, $\Rt(A)\dotplus \Rt(T)^{\perp}=P_T(\Rt(A)) \oplus \Rt(T)^{\perp}=\Rt(T)\oplus \Rt(T)^{\perp}$ is dense in $\HH.$ 
\end{dem}

Given $T \in \PL,$ define the set
$$\Phi(T):=\{ E=P_{\Rt(A) \pl \Rt(T)^{\perp}} \mbox{ such that } A \mbox{ is optimal for } T \}.$$

By Proposition \ref{propMP2}, every $E \in \Phi(T)$ is a densely defined semiclosed projection with closed nullspace. Moreover, $\DD(E)= \Rt(T)\oplus \Rt(T)^{\perp} $ and $\ol{\Rt(E)}=\ol{\Rt(T^*)}.$ Also, there exists a unique $E \in \Phi(T)$ if and only if  $\ol{\Rt(T^*)} \cap \N(T^*)=\{0\},$ see \cite[Proposition 4.1]{Arias}.

\section{$B$-symmetric projections}
A densely defined operator $T$ is \emph{symmetric} if $T \subset T^*$ and it is \emph{selfadjoint} if $T=T^*,$ i.e., $T$ is symmetric and $\DD(T)=\DD(T^*).$

\begin{Def} Let $B \in L(\HH)^s$ and $E$ be a densely defined projection. We say that $E$ is \emph{$B$-symmetric} if $BE$ is symmetric and it is \emph{$B$-selfadjoint} if $BE$ is selfadjoint.
\end{Def}

Since $B$ is bounded, $\DD(BE)=\DD(E)$ and $(BE)^*=E^*B.$ Therefore $E$ is $B$-symmetric if and only if $BEx=E^*Bx$ for every $x \in \DD(E)$ and $E$ is $B$-selfadjoint if and only if $BE=E^*B.$

\begin{prop} \label{EBself} Let $B \in L(\HH)^s$ and $E$ be a densely defined projection in $\HH.$ 
Then $E$ is $B$-symmetric if and only if $\N(E) \subseteq (B\Rt(E))^{\perp}.$
\end{prop}

\begin{dem}  Suppose that $E$ is $B$-symmetric. Let $x \in \N(E)$ and $y \in \DD(E).$ Then $$\PI{x}{BEy}=\PI{x}{E^*By}=\PI{Ex}{By}=0$$ so, $\N(E) \subseteq  (B\Rt(E))^{\perp}.$ 
Conversely, suppose that $\N(E) \subseteq (B\Rt(E))^{\perp}=B^{-1}(\N(E^*)).$ Then $B\N(E) \subseteq \N(E^*)$ and $B\Rt(E)\subseteq \ol{B\Rt(E)} \subseteq \N(E)^{\perp}=\Rt(E^*).$ Then $B\DD(E)=B\Rt(E)+B\N(E)\subseteq \Rt(E^*) \dotplus \N(E^*)= \DD(E^*).$ 
Then $E^*Bx=0$ for every $x \in \N(E)$ and  $E^*Bx=Bx=BEx$ for every $x \in \Rt(E).$ Finally, if $x \in \DD(E)$ then $Bx, BEx, B(I-E)x \in \DD(E^*)$ and  $E^*Bx=E^*BEx+E^*B(I-E)x=BEx.$ Then $E$ is $B$-symmetric.
\end{dem}

\begin{prop} \label{propCuasi} Let $B \in L(\HH)^s$ and $\St$ a (not necessarily closed) subspace. 
There exists a $B$-symmetric projection onto $\St$ if and only if  $$\HH=\ol{\St+(B\St)^{\perp}}.$$
In this case, $\St \cap (B\St)^{\perp}= \St \cap \N(B).$
\end{prop}
\begin{dem}
	Suppose that $E$ is a $B$-symmetric projection onto $\St.$ Then, by Proposition \ref{EBself}, $\N(E) \subseteq (B\St)^{\perp}$ and $\DD(E)=\Rt(E) \dotplus \N(E) \subseteq \St + (B\St)^{\perp}.$ Therefore $\HH=\ol{\DD(E)}\subseteq \ol{\St + (B\St)^{\perp}}.$ 
	
	Conversely, suppose that $\HH=\ol{\St+(B\St)^{\perp}}=\ol{\ol{\St}+(B\St)^{\perp}}.$ Let $\mc{L'}:=\ol{\St} \cap (B\St)^{\perp},$ then $\ol{\St}+(B\St)^{\perp}=\ol{\St} \dotplus (B\St)^{\perp} \cap\mc{L'}^{\perp},$ hence $\St \cap ((B\St)^{\perp} \cap \mc{L'}^{\perp})=\{0\}.$ 
	
	Define $E:=P_{\St {\mathbin{\!/\mkern-3mu/\!}} (B\St)^{\perp} \cap \mc{L'}^{\perp}}.$ Then $\ol{\DD(E)}=\ol{ \St \dotplus ((B\St)^{\perp} \cap \mc{L'}^{\perp})}=\ol{\ol{\St} \dotplus (B\St)^{\perp}\cap\mc{L'}^{\perp}}=\ol{\ol{\St}+(B\St)^{\perp}}=\HH.$
	Therefore $E$ is a densely defined projection with (closed) nullspace contained in $(B\St)^{\perp}$  and, by Proposition \ref{EBself}, $E$ is $B$-symmetric. 
	
	Finally, the inclusion $\St \cap \N(B) \subseteq \St \cap (B\St)^{\perp}$ always holds. On the other hand, let $x \in \St \cap (B\St)^{\perp}$ and $y \in \DD(E)$ then 
	$$\PI{Bx}{y}=\PI{BEx}{y}=\PI{x}{BEy}=0=\PI{Bx}{Ey}.$$ Then $Bx \in \DD(E)^{\perp}=\{0\}$ whence $x \in \St \cap \N(B).$
\end{dem}

When the projection $E$ is semiclosed, the $B$-symmetry can be given in terms of bounded operators.

\begin{prop} \label{EBselfSC} Let $B \in L(\HH)^s$ and $E\in \Q_{SC}.$ Then $E$ is $B$-symmetric if and only if $P_\G$ commutes with $\G B \G,$ where $\G \in L(\HH)^+$ and $P_{\G} \in \mc{P}$ are as in Corollary \ref{ddprojection}. 
\end{prop}

\begin{dem}  Suppose that $E$ is  $B$-symmetric, then $B\G P_\G x=BE\G x=E^*B\G x$ for every $x \in \HH.$ Then  $\G B\G P_\G x=\G E^*B \G x$ for every $x \in \HH.$ Now, since $E\G = \G P_\G,$ it follows that $\G E^* \subset (E\G)^* = P_\G^* \G.$ Therefore
$\G B\G P_\G x=P_\G^*\G B \G x$ for every $x \in \HH,$ i.e., $P_\G$ is $\G B \G$-selfadjoint and since $P_\G$ is selfadjoint, then $P_\G$ commutes with $\G B \G.$
	
Conversely,  if $P_\G$ commutes with $\G B \G,$ $\N(P_\G) \subseteq (\G B \G)^{-1}(\N(P_\G)).$ Therefore $$\N(E)=\G(\N(P_\G)) \subseteq  B^{-1}(\Rt(\G P_\G)^{\perp}) \subseteq B^{-1}(\Rt(E)^{\perp}).$$ Then, by Proposition \ref{EBself}, $E$ is $B$-symmetric.
\end{dem}

\section{Quasi and weak complementability}

The complementability of an operator $B \in L(\HH)$ with respect to two given closed subspaces $\St$ and $\T$ of $\HH$ was studied for matrices by Ando \cite{AndoSchur} and Carlson and Haynsworth \cite{Carlson}. These ideas were extended to operators in Hilbert spaces in \cite{CMSSzeged,AntCorSto06}.

\bigskip
\begin{Def} Let $B \in L(\HH)^s$ and $\St \subseteq \HH$ be a closed subspace. Then $B$ is $\St$\emph{-complementable} if
	$$\HH=\St+ (B\St)^{\perp}.$$
\end{Def}

In \cite{CMSSzeged} it was shown that $B$ is $\St$-complementable if and only if there exists a $B$-selfadjoint projection onto $\St;$ i.e., the set
$$\mc{P}(B,\St):=\{ Q \in \Q: \Rt(Q)=\St, \ BQ \in L(\HH)^s \}$$ is not empty.

\begin{prop} [\cite{CMSSzeged}] Let $B \in L(\HH)^s$ and $\St \subseteq \HH$ be a closed subspace. If the matrix decomposition of $B$ according to the orthogonal decomposition $\HH=\St \oplus \St^{\perp}$ is given by  
	\begin{equation} \label{matrixW}
	B=\begin{bmatrix}
	a  & b \\ 
	b^* & c \\
	\end{bmatrix},
	\end{equation}
	then $B$ is $\St$-complementable if and only if $\Rt(b) \subseteq \Rt(a).$ 
\end{prop}

One way of generalizing the concept of complementability is to consider $B$-symmetric closed projections onto $\St.$
\begin{Def} Let $B \in L(\HH)^s$ and $\St \subseteq \HH$ be a closed subspace. We say that the pair $(B,\St)$ is \emph{quasi-complementable} if there exists a $B$-symmetric closed projection onto $\St.$
\end{Def}

The set of quasi-complementable pairs was studied in \cite{Cuasi} for a positive weight $B.$ Many results stated in \cite{Cuasi} hold also in the selfadjoint case.

\begin{prop} \label{EBselfC4} Let $B \in L(\HH)^s$ and $\St \subseteq \HH$ be a closed subspace. The pair $(B,\St)$ is \! quasi-complementable if and only if $\ol{B\St} \cap \St^{\perp}=\{0\}.$
\end{prop}

\begin{dem} It follows from the definition of quasi-complementability and Proposition \ref{propCuasi}.
\end{dem}

\begin{prop} \label{propaxb} Let $B \in L(\HH)^s,$ $\St \subseteq \HH$ be a closed subspace and $E$ be a densely defined closed projection onto $\St.$ Then $E$ is $B$-symmetric if and only if $$ax \subset b,$$
	where $a, b$ are as in \eqref{matrixW} and $x$ is as in \eqref{matrixE}. 
\end{prop}

\begin{dem} The result follows by using arguments similar to those found in \cite[Proposition 2.2]{Cuasi}.
\end{dem}

\begin{cor}  Let  $\St \subseteq \HH$ be a closed subspace  and $B \in L(\HH)^s$  with representation as in \eqref{matrixW}. Then $(B,\St)$ is quasi-complementable if and only if the equation 
\begin{equation} \label{cuasieq}
b^*=xa
\end{equation} admits a densely defined closed solution $x_0: \DD(x_0) \subseteq \St \ra \St^{\perp}$ with densely defined adjoint.
\end{cor}

\begin{dem} If $(B,\St)$ is quasi-complementable then there exists $E$ a $B$-symmetric closed projection onto $\St.$ Suppose that the matrix decomposition of $E$ is as in \eqref{matrixE}. Then, by Proposition \ref{propaxb}, $ax \subset b.$ Then $b^* \subset (ax)^* = x^*a.$ Therefore, $b^*=x^*a$ and $x^*$ is a densely defined closed  solution of \eqref{cuasieq} with densely defined adjoint.
	
Conversely, let $x_0: \DD(x_0) \subseteq \St \ra \St^{\perp}$ be a densely defined closed solution of \eqref{cuasieq}  with densely defined adjoint. Then $b=(x_0a)^* \supset ax_0^*.$ Set $E:=\begin{bmatrix}
	I & x_0^* \\ 
	0 & 0 \\
\end{bmatrix},$ then $E$ is a densely defined closed projection with range $\St$ and, by Proposition \ref{propaxb}, $E$ is $B$-symmetric. Hence, the pair $(B,\St)$ is quasi-complementable.
\end{dem}

\bigskip
 A different way of extending the concept of complementability was given in \cite{AntCorSto06}, where Antezana et al. defined the notion of \emph{weak complementability} to study the Schur complement in this context. We use these ideas when $\St=\T$ and $B \in L(\HH)^s.$ 
\begin{Def} Let  $\St \subseteq \HH$ be a closed subspace, and $B \in L(\HH)^s$  with representation as in \eqref{matrixW}. Then $B$ is $\St$\emph{-weakly complementable} if
	$$\Rt(b) \subseteq \Rt(\vert a \vert^{1/2}).$$
\end{Def}
The notion of weak complementability is distinct to the notion of complementability only in the infinite dimensional setting.
Every positive operator $B$ is $\St$-weakly complementable, and if $B$ is $\St$-weakly complementable for every closed subspace $\St \subseteq \HH$ then $B$ is semidefinite, see \cite[Proposition 3.1]{Contino4}.

The next proposition gives an operator characterization of the $\St$-weak complementability as the solution of a Riccati type equation \cite{Riccati}.
\begin{prop} \label{PropWC2} Let $B \in L(\HH)^s$ and $\St$ be a closed subspace of $\HH.$ Then 
	$B$ is $\St$-weakly complementable if and only if there exists a positive solution of the equation $$BP_{\St}B=XP_{\St}X^*.$$
\end{prop}

\begin{dem} Suppose that the matrix decomposition of $B$ induced by $\St$ is as in \eqref{matrixW} and $a=u\vert a \vert$ is the polar decomposition of $a.$ Let $f \in L(\St^{\perp},\St)$ be the reduced solution of $b=\vert a \vert^{1/2}x.$ 
	
	Let $A:= \begin{bmatrix} 
	\vert a\vert & u\vert a \vert^{1/2} f \\ 
	f^* \vert a \vert^{1/2} u & f^*f \\  
	\end{bmatrix}.$ Then, by Lemma \ref{LemmaPositive}, $A \geq 0,$ because $uf$ is the reduced solution of $\vert a \vert^{1/2}x = u \vert a \vert^{1/2} f$ and $f^*f=f^*uuf.$  It easily follows that $BP_{\St}B=AP_{\St}A.$ 
	
	Conversely, suppose that $BP_{\St}B=AP_{\St}A$ with $A \geq 0$ and $A=\begin{bmatrix} 
	a_{11} & a_{12} \\ 
	a_{12}^* & a_{22} \\  
	\end{bmatrix}.$ Then, 
	\begin{eqnarray}
	a^2=a_{11}^2  \label{equation1} \\
	ab=a_{11}a_{12} \label{equation2} \\
	b^*b=a_{12}^*a_{12}. \label{equation3}
	\end{eqnarray}
	From \eqref{equation1}, $a_{11}=\vert a \vert.$ Thus, by  \eqref{equation2}, $ua_{11}b=ab=a_{11}a_{12}.$  Then $a_{11}b=ua_{11}a_{12}=a_{11}ua_{12}.$ So that $$a_{11}^{1/2}(b-ua_{12})=0=(b-ua_{12})^*a_{11}^{1/2}.$$
	Then, from \eqref{equation3}, $a_{12}^*a_{12}=b^*b=(b-ua_{12}+ua_{12})^*(b-ua_{12}+ua_{12})=(b-ua_{12})^*(b-ua_{12})+a_{12}^*a_{12}.$  The last equality follows from  Lemma \ref{LemmaPositive}; in fact, if $a_{12}=a_{11}^{1/2}g,$ then $ua_{12}=ua_{11}^{1/2}g=a_{11}^{1/2}ug,$ because $u$ and $a_{11}^{1/2}$ commute. Then $(b-ua_{12})^*ua_{12}=(b-ua_{12})^*a_{11}^{1/2}ug=0.$ Therefore,
	$$b^*b=(b-ua_{12})^*(b-ua_{12})+b^*b.$$ So that
	 $\vert b-ua_{12} \vert =0$ or $b=ua_{12}=a_{11}^{1/2}ug$ and $\Rt(b) \subseteq \Rt(a_{11}^{1/2})=\Rt(\vert a \vert^{1/2}).$
\end{dem}

\bigskip
Let $\St$ be a closed subspace of $\HH$ and let $E$ be a densely defined projection with $\N(E)=\St^{\perp}.$ Then, by Proposition \ref{proprepres}, the matrix representation of $E$ according to the orthogonal decomposition $\HH=\St \oplus \St^{\perp}$ is
\begin{equation} \label{ESCDecomp}
E=\begin{bmatrix}
I & 0 \\ 
y & 0 \\
\end{bmatrix} \mbox {  where } \DD(E)=\DD(y) \oplus \St^{\perp},
\end{equation}
with $\DD(y)=P_{\St}(\Rt(E))$ dense in $\St$ and $y=P_{\St^{\perp}}(P_{\St}P_{\ol{\Rt(E)}}|_{\Rt(E)})^{-1}.$ 

To characterize the $\St$-weak complementability of $B$ in terms of projections we need the following lemma. 

\begin{lema} \label{ProjectionBself} Let $B \in L(\HH)^s,$ $\St$ be a closed subspace of $\HH$ and $E$ be a densely defined projection with $\N(E)=\St^{\perp}.$ Suppose that the matrix decomposition of $B$ and $E$ are as in \eqref{matrixW} and \eqref{ESCDecomp}, respectively. 
Then $EB \in L(\HH)^s$ if and only if $ya=b^*$ and $yb \in L(\St^{\perp})^s.$ 
Moreover,  $B\St \subseteq \Rt(E)$ if and only if $ya=b^*.$ 
\end{lema}

\begin{dem} Suppose that $EB \in L(\HH)^s.$ Then $EB=\begin{bmatrix}
	a & b \\ 
	ya & yb \\
	\end{bmatrix}=(EB)^*=\begin{bmatrix}
	a & (ya)^* \\ 
	b^* & (yb)^* \\
	\end{bmatrix}.$ Therefore $ya=b^*$ and $yb \in L(\St^{\perp})^s.$ The converse follows in a similar way.
	
	Moreover, if $B\St \subseteq \Rt(E)$ then 
$$0=(I-E)Bs=\begin{bmatrix}
	0 & 0 \\ 
	-y & I \\
\end{bmatrix}\begin{bmatrix}
a & b \\ 
b^* & c \\
\end{bmatrix}\m{s}{0}=-yas+b^*s \mbox{ for every } s \in \St.$$ Then $b^*=ya.$

Conversely, if $b^*=ya$ then $B\St \subseteq \Rt(a)+\Rt(b^*) \subseteq \DD(y) \oplus \St^{\perp}=\DD(E)=\DD(I-E)$ and $(I-E)Bs=-yas+b^*s=0 \mbox{ for every } s \in \St.$ Hence, $B\St \subseteq \Rt(E).$
\end{dem}

\begin{cor}  Let $B \in L(\HH)^s,$ $\St$ be a closed subspace of $\HH$ and $E$ be a densely defined projection with $\N(E)=\St^{\perp}.$ If  $EB \in L(\HH)^s$ then $B\St\cap \St^{\perp}=\{0\}.$
\end{cor}

\begin{cor} \label{corEBself} Let $B \in L(\HH)^s,$ $\St$ be a closed subspace of $\HH$ and $E\in \Q_{SC}$ with $\N(E)=\St^{\perp}.$ Suppose that the matrix decomposition of $B$ and $E$ are as in \eqref{matrixW} and \eqref{ESCDecomp}, respectively. Then $EB \in L(\HH)^s$ and $\Rt(\vert a \vert^{1/2}) \subseteq \DD(E)$
	if and only if $ya=b^*$ and $\Rt(\vert a \vert^{1/2})\subseteq \DD(y).$  
\end{cor}

\begin{dem} Suppose that  $EB \in L(\HH)^s$ and $\Rt(\vert a \vert^{1/2}) \subseteq \DD(E).$ Then, by Lemma \ref{ProjectionBself}, $ya=b^*.$ Also, since $\DD(E)=\DD(y) \oplus \St^{\perp},$ it follows that $\Rt(\vert a \vert^{1/2}) \subseteq \DD(y).$ 
	
	Conversely, suppose that $ya=b^*$ and $\Rt(\vert a \vert^{1/2})\subseteq \DD(y).$  Then  $\Rt(\vert a \vert^{1/2}) \subseteq \DD(y) \oplus \St^{\perp}=\DD(E).$ By Corollary \ref{corKaufman2}, $E\vert a \vert^{1/2}\in L(\St, \HH),$ so that $y\vert a \vert^{1/2} \in L(\St, \St^{\perp}).$ Then, if $a=u\vert a \vert$ is the polar decomposition of $a$, 
	$$yb=y(ya)^*=y (y\vert a \vert^{1/2}u\vert a \vert^{1/2})^*=y\vert a \vert^{1/2}u(y\vert a \vert^{1/2})^*.$$ Hence $yb \in L(\St^{\perp})^s$ and, by Lemma \ref{ProjectionBself}, $EB \in L(\HH)^s.$
\end{dem}

The next theorem characterizes the $\St$-weak complementability of a selfadjoint operator in terms of semiclosed projections. See also \cite[Theorem 3.14]{Contino4}.

\begin{thm} \label{SWC} Let $B \in L(\HH)^s$ and $\St$ be a closed subspace of $\HH.$ Suppose that the matrix decomposition of $B$ is as in \eqref{matrixW}.  Then $B$ is $\St$-weakly complementable if and only if there exists $E \in \Q_{SC}$ with $\N(E)=\St^{\perp}$ such that $EB \in L(\HH)^s$ and $\Rt(\vert a \vert^{1/2}) \subseteq \DD(E).$
\end{thm}

\begin{dem} Suppose that $B$ is $\St$-weakly complementable. If the matrix decomposition of $B$ induced by $\St$ is as in \eqref{matrixW}, let $f$ be the reduced solution of $b=\vert a\vert^{1/2} x$ and $a=u \vert a \vert$ be the polar decomposition of $a.$
Write $(\vert a \vert^{1/2})^{\dagger}$ for the Moore-Penrose inverse of $\vert a \vert^{1/2}$ and set $$E:=\begin{bmatrix} 
I & 0 \\ 
f^*u(\vert a \vert^{1/2})^{\dagger} & 0\\  
\end{bmatrix}.$$
	
Then $\DD(E)=\DD((\vert a \vert^{1/2})^{\dagger}) \oplus \St^{\perp}$ is a semiclosed subspace of $\HH$ (because it is the sum of two semiclosed subspaces), $E$ is a densely defined projection with $\N(E)=\St^{\perp}$ and $\Rt(\vert a \vert^{1/2}) \subseteq \DD(E).$ 
On the other hand, since $\Rt(B) \subseteq \Rt(\vert a \vert^ {1/2}) \oplus \St^{\perp},$ the product $(I-E)B$ is well defined. Moreover,	
	\begin{align*}
	(I-E)B&=\begin{bmatrix} 
	0 & 0 \\ 
	-f^*u(\vert a \vert^{1/2})^{\dagger} & I\\  
	\end{bmatrix} \begin{bmatrix} 
	\vert a\vert^{1/2} u \vert a\vert^{1/2} & \vert a\vert^{1/2}f \\ 
	f^*\vert a\vert^{1/2} &c\\  
	\end{bmatrix} \\
	& = \begin{bmatrix} 
	0 & 0 \\ 
	0 &c-f^*uf\\  
	\end{bmatrix} \in L(\HH)^s.
	\end{align*}
	So that $EB \in L(\HH)^s.$

	Finally, if $\Gamma:=\begin{bmatrix} 
		\vert a \vert^{1/2}+P_{\N(a)} & 0 \\ 
		0 & I\\  
	\end{bmatrix},$ then  $\Gamma\in L(\HH)^{+},$  $\Rt(\Gamma)=\DD(E)$ and

	$$E \Gamma= \begin{bmatrix} 
	\vert a \vert^{1/2}+P_{\N(a)} & 0 \\ 
	f^*u & 0\\  
	\end{bmatrix} \in L(\HH).$$ By Corollary \ref{corKaufman}, $E$ is semiclosed.

	Conversely, suppose that there exists $E\in \Q_{SC}$ with $\N(E)=\St^{\perp}$ such that $EB \in L(\HH)^s$ and $\Rt(\vert a \vert^{1/2}) \subseteq \DD(E).$ Suppose that the matrix decomposition of $E$ is that as in \eqref{ESCDecomp}. By Corollary \ref{corKaufman2}, $E\vert a \vert^{1/2} \in L(\St, \HH),$ and then $y\vert a \vert^{1/2} \in L(\St, \St^{\perp}).$  By Lemma \ref{ProjectionBself}, $ya=b^*$ and $yb \in L(\St^{\perp})^s.$ Since $ya=b^*,$ we also have that $y \vert a \vert^{1/2} u \vert a \vert^{1/2}= b^*.$ Then $b=\vert a \vert^{1/2} u (y \vert a \vert^{1/2})^*,$ $\Rt(b) \subseteq \Rt(\vert a \vert^{1/2})$ and $B$ is $\St$-weakly complementable. 
\end{dem}

Suppose that $B$ is $\St$-weakly complementable. If the matrix decomposition of $B$ induced by $\St$ is as in \eqref{matrixW}, let $f$ be the reduced solution of $b=\vert a\vert^{1/2} x$ and $a=u \vert a \vert$ be the polar decomposition of $a.$ Set
\begin{equation} \label{E0}
E_0:=\begin{bmatrix} 
I & 0 \\ 
y_0 & 0\\  
\end{bmatrix},
\end{equation}
with $y_0:=f^*u(\vert a \vert^{1/2})^{\dagger}.$ Then, by the proof of Theorem \ref{SWC}, $E_0 \in \Q_{SC}$ with $\N(E_0)=\St^{\perp},$ $E_0B \in L(\HH)^s$ and $\DD(E_0)=\DD((\vert a \vert^{1/2})^{\dagger}) \oplus \St^{\perp}.$ Define 
 $$\mc{P}^*(B,\St):=\{ E \in  \Q_{SC} :  \ \N(E)=\St^{\perp}, \ EB\in L(\HH)^s \mbox{ and } \Rt(\vert a \vert^{1/2}) \subseteq \DD(E)\}.$$ Let $E \in \mc{P}^*(B,\St)$ with matrix decomposition as in \eqref{ESCDecomp}. Then, it can be proved that
 $$\Rt(y_0) \subseteq \Rt(y).$$
 Furthermore, $y\vert a \vert^{1/2}=y_0\vert a \vert^{1/2}.$
 
 In what follows, we characterize the subset of projections with fixed domain $\DD(E_0).$ More precisely, consider
\begin{align*}
\mc{P}_{0}^*(B,\St)&:=\{E \in \mc{P}^*(B,\St) :\DD(E)=\DD(E_0)\}.
\end{align*}
From  Lemma \ref{ProjectionBself} and Corollary \ref{corEBself},
$\mc{P}_{0}^*(B,\St)=\{ E \in \Q_{SC} :  \DD(E)= \DD(E_0), \ \N(E)=\St^{\perp} \mbox{ and } B\St \subseteq \Rt(E)\}.$

Clearly, $\mc{P}_{0}^*(B,\St) \subseteq \mc{P}^*(B,\St)$ and $\mc{P}^*_{0}(B,\St)$ is not empty because $E_0 \in \mc{P}^*_{0}(B,\St).$

\begin{thm} Let $B \in L(\HH)^s$ and $\St$ be a closed subspace of $\HH$ such that $B$ is $\St$-weakly complementable. Then
	$$\mc{P}_{0}^*(B,\St)=E_0+\{W\in SC(\HH): \DD(E_0) \subseteq \DD(W), \ \Rt(W) \subseteq \St^{\perp} \mbox{ and } B\St\dotplus\St^{\perp} \subseteq \N(W)\}.$$
\end{thm}

\begin{dem} Let $E=E_0+W,$ with $W\in SC(\HH)$ such that  $\DD(E_0) \subseteq \DD(W), \ \Rt(W) \subseteq \St^{\perp} \mbox{ and } B\St\dotplus\St^{\perp} \subseteq \N(W).$ Then $E_0+W \in SC(\HH)$ and $\DD(E_0+W)=\DD(E_0) \cap \DD(W) = \DD(E_0).$ Observe that, $\Rt(E_0+W) \subseteq \Rt(E_0) \dotplus \St^{\perp}=\DD(E_0)=\DD(E_0+W).$ Also, since $\Rt(W) \subseteq \N(W) \cap \N(E_0),$ $W^2=0,$ $E_0W=0$ and, from $\Rt(I-E_0) \subseteq \N(W),$ $WE_0=W.$ Hence $E_0+W \in \Q_{SC}.$ Furthermore, $\N(E_0+W)=\St^{\perp}.$  In fact, it is clear that $\St^{\perp} \subseteq \N(E_0+W)$ and if $h \in \N(E_0+W)\subseteq \DD(E_0+W)=\DD(E_0)$ then $E_0h=-Wh,$ so that $E_0h \in \Rt(E_0) \cap \Rt(W) \subseteq  \Rt(E_0) \cap \St^{\perp}=\{0\}.$ Therefore $h \in \N(E_0)=\St^{\perp}$ and $\N(E_0+W)\subseteq \St^{\perp}.$ 
Also, $\Rt(B) \subseteq \DD(E_0) = \DD(E_0+W)$ and $B\St \subseteq \Rt(E_0+W).$ In fact, if $h \in B\St,$ since $E_0B \in L(\HH)^s,$ by Lemma \ref{ProjectionBself}, $h \in \Rt(E_0) \subseteq \DD(E_0)= \DD(W+E_0).$ Hence, $(E_0+W)h=E_0h+Wh=h,$ because $h \in \N(W).$ Then  	$E_0+W \in \mc{P}_{0}^*(B,\St).$

Conversely, let $E \in 	\mc{P}_{0}^*(B,\St)$ and define $W:=E-E_0.$ Then $W \in SC(\HH)$ and $\DD(W)=\DD(E) \cap \DD(E_0) = \DD(E_0).$ Also, $\Rt(W)=\Rt((I-E_0)+(E-I)) \subseteq \St^{\perp}.$ It is clear that $\St^{\perp} \subseteq \N(W)$ and, if $h \in B\St,$ by Lemma \ref{ProjectionBself},  $h  \in \Rt(E) \cap \Rt(E_0) \subseteq \DD(E_0)$. Then $Wh=Eh-E_0h=h-h=0.$ Therefore $B\St \dotplus \St^{\perp} \subseteq \N(W)$ and $E=E_0+W.$ 
\end{dem}

\subsection{Comparison between the notions of quasi and weak complementability}

The next examples show that  quasi-complementability does not imply weak complementability and viceversa.

\begin{example} \label{ex1} Let $\St \subseteq \HH$ be a closed subspace such that $\dim(\St)=\dim(\St^{\perp})=\infty.$ Take $a \in L(\St)^+$ such that $\Rt(a)$ is not closed and $\ol{\Rt(a)}=\St,$ and take $b: \St^{\perp} \ra \St$ such that $b$ is invertible. 
Consider $B=\begin{bmatrix}
a & b \\ 
b^* & c \\
\end{bmatrix}$ for any $c \in L(\St^{\perp})^s.$ Then $B$ is not $\St$-weakly complementable because $\Rt(b)=\St \not \subseteq \Rt(a^{1/2}).$ 

Now, we want to show that $\ol{B\St} \cap \St^{\perp}=\{0\}.$ Let $h \in \ol{B\St} \cap \St^{\perp}.$ Then there exists $\{s_n\}_{n \geq 1} \subseteq \St$ such that $$h=\underset{n \ra \infty}{\lim}Bs_n=\underset{n \ra \infty}{\lim} (as_n+b^*s_n)=\underset{n \ra \infty}{\lim} P_{\St^{\perp}}(as_n+b^*s_n)=\underset{n \ra \infty}{\lim} b^*s_n.$$ Then $(b^*)^{-1}h=\underset{n \ra \infty}{\lim}s_n$ and $0=\underset{n \ra \infty}{\lim}as_n.$ Therefore $a(b^*)^{-1}h=0.$ Since $a$ is injective it follows that $(b^*)^{-1}h=0$ then $h=0.$ Hence, by Proposition \ref{EBselfC4}, the pair $(B,\St)$ is quasi-complementable.
\end{example}

\begin{example} See \cite[Example 2.14]{Cuasi}. Let $B \in L(\HH)^+$ be such that $\Rt(B)$ is not closed. Let $x \in \ol{\Rt(B)} \setminus \Rt(B)$ and define the closed subspace $\St$ such that $\St^{\perp}=span \ \{ x \}.$ Clearly, $B$ is $\St$-weakly complementable. On the other hand, $(B\St)^{\perp}=B^{-1}(\St^{\perp})=B^{-1}(span\ \{ x \})=B^{-1}(span\  \{ x \}  \cap \Rt(B))=B^{-1}(\{0\})=\N(B).$ Then $\ol{BS}=\ol{\Rt(B)}$ and $\St^{\perp} \cap \ol{B\St} \not = \{0\}.$ Therefore, by Proposition \ref{EBselfC4}, the pair $(B, \St)$ is not quasi complementable.
\end{example}

Let $B \in L(\HH)^s$ and let $\St$ be a closed subspace of $\HH$ and suppose that $\St = \St_{+} \ \oplus_{B} \ \St_{-}$ is any decomposition as in \eqref{Sdecomp}. In \cite[Proposition 3.2]{Contino4} it was shown that
$B$ is $\St$-weakly complementable if and only if there exist $B_1, B_2, B_3 \in L(\HH)^s,$ such that
\begin{equation} \label{eqWC}
B=B_1+B_2-B_3,
\end{equation} 
where $B_2, B_3 \geq 0,$ $\St \subseteq \N(B_1),$ $\St_- \subseteq \N(B_2),$ $\St_+ \subseteq \N(B_3).$

In the following proposition we characterize the $B$-symmetric closed projections onto $\St$ for $B \in L(\HH)^s$ when $B$ is $\St$-weakly complementable. Some of these results where stated in \cite{Cuasi} for a positive weight $B.$ To extend these results to the selfadjoint case the notion of semiclosed projections turns out to be useful.

\begin{prop} \label{EBselfC2} Let $B \in L(\HH)^s,$ $E$ be a densely defined closed projection of $\HH$ onto $\St$ and suppose that $B$ is $\St$-weakly complementable.  If $E$ is $B$-symmetric, then $BE$ admits a bounded selfadjoint extension to $\HH.$ 
Moreover, if $\St=\St_{1} \ \oplus_{B} \ \St_{2}$ and  $B=B_1+B_2-B_3,$ are any decompositions as in \eqref{Sdecomp} and  \eqref{eqWC}, respectively, then $$\ol{BE}=(BE)^*=B_2^{1/2}P_{\M_2}B_2^{1/2}-B_3^{1/2}P_{\M_3}B_3^{1/2} \in L(\HH)^s,$$ where $\M_2=\ol{B_2^{1/2}(\St_+)}$ and $\M_3=\ol{B_3^{1/2}(\St_-)}.$ 
\end{prop}

To prove this proposition we need the following lemma. 

\begin{lema} \label{EBselfC} Let $B \in L(\HH)^s$ and $E$ be a densely defined closed projection of $\HH$ onto $\St.$ Suppose that $B$ is $\St$-weakly complementable, $\St=\St_{+} \ \oplus_{B} \ \St_{-}$ is any decomposition as in \eqref{Sdecomp} and $B=B_1+B_2-B_3$ is any decomposition as in \eqref{eqWC}.  If $E$ is $B$-symmetric, then $E$ admits a factorization $E=E_++E_-,$ where $E_+$ and $E_-$ are semiclosed projections, with $\DD(E)=\DD(E_+)=\DD(E_-),$ $\Rt(E_+)=\St_+,$ $\Rt(E_-)=\St_-,$  $E_-E_+=E_+E-=0,$ $E_+$ is $B_2$-symmetric and $E_-$ is  $B_3$-symmetric.
\end{lema}

\begin{dem} 
Since $E$ is $B$-symmetric, by Proposition  \ref{EBself}, $E=P_{\St {\mathbin{\!/\mkern-3mu/\!}} \T},$ with $\T$ a closed subspace such that $\T \subseteq B^{-1}(\St^{\perp})$ and $\DD(E)= \St \dotplus \T =\St_+ \dotplus \St_- \dotplus \T.$
Let $E_+:=P_{\St_+ {\mathbin{\!/\mkern-3mu/\!}} \T + S_-}$ and $E_-:=P_{\St_- {\mathbin{\!/\mkern-3mu/\!}} \T + S_+}.$ Then $E_\pm=P_{\St_\pm}E$ and $E_\pm$ are $B$-symmetric.  Let us prove these assertions for $E_+;$ the other case is similar. First observe that $\DD(P_{\St_+}E)=\DD(E)$ and $\Rt(P_{\St_+}E)\subseteq \St_+ \subseteq \DD(E).$ Then $(P_{\St_+}E)^2=P_{\St_+}EP_{\St_+}E=P_{\St_+}E$ so that $P_{\St_+}E$ is a densely defined projection. On the other hand, $\Rt(P_{\St_+}E)=P_{\St_+}\Rt(E)=\St_+$ and $\N(P_{\St_+}E)=\N(E)+\St \cap \St_+^{\perp}=\T + \St_-$ (see Lemma \ref{Sdecomp}). Therefore $E_+=P_{\St_+}E$ and, since $\N(E_+)$ and $\Rt(E_+)$ are semiclosed subspaces, $E_+$ is a semiclosed projection. Also, $\DD(E)=\DD(E_+)=\DD(E_-)$ and $E=P_{\St}E=P_{\St_+}E+P_{\St_-}E=E_++E_-$ with $E_+E_-=E_-E_+=0.$

Now, let us see that $E_+$ is $B$-symmetric. Note that $P_\St E=E,$ so that $E^*P_\St=E^*.$ Let $x \in \DD(E_+)=\DD(E),$ since $E_+x \in \DD(E),$ 
\begin{align*}
BE_+x&=BEE_+x=E^*BE_+x=E^*BP_{\St_+}Ex=E^*P_{\St}BP_{\St_+}Ex=E^*P_{\St}BP_{\St}P_{\St_+}Ex\\
&=E^*G_{B,\St}P_{\St_+}Ex=E^*P_{\St_+}G_{B,\St}Ex=E^*P_{\St_+} P_\St B Ex=E^*P_{\St_+} B Ex\\
&=E^*P_{\St_+} E^*Bx=(P_{\St_+}E)^*E^*Bx=E_+^*E^*Bx=E_+^*Bx,
\end{align*}
where we used  that $G_{B,\St}P_{\St_+}=P_{\St_+}G_{B,\St}$ and $E_+^*=(EE_+)^* \supset E_+^*E^*.$
Hence $BE_+ \subset E_+^*B.$ 

Note that $BE_+=B_1E_++B_2E_+-B_3E_+=B_2E_+,$ because $\Rt(E_+) \subseteq \St =\N(B_1)$ and $\Rt(E_+)=\St_+ \subseteq \N(B_3).$ Then
$$B_2E_+=BE_+ \subset E_+^*B = (BE_+)^*=(B_2E_+)^*=E_+^*B_2.$$ Hence $E_+$ is $B_2$-symmetric.
\end{dem}

{\vspace{1ex}\noindent{\it Proof of Proposition \ref{EBselfC2}.}\hspace{0.5em}} 
By Lemma \ref{EBselfC}, $E$ admits a factorization in the form $E=E_++E_-,$ where $E_+$ is a $B_2$-symmetric semiclosed projection with range $\St_+$ and $E_-$ is a $B_3$-symmetric semiclosed projection with range $\St_-.$ 

Let us show that $P_{\M_2}B_2^{1/2}=B_2^{1/2}E_+$ in $\DD(E_+).$ In fact, if $x \in \DD(E_+)$ then $$P_{\M_2}B_2^{1/2}x=P_{\M_2}B_2^{1/2}E_+x+P_{\M_2}B_2^{1/2}(I-E_+)x.$$ Since $E_+x \in \St_+,$ $P_{\M_2}B_2^{1/2}E_+x=B_2^{1/2}E_+x.$ Also, $P_{\M_2}B_2^{1/2}(I-E_+)x=0$ because $B_2^{1/2}(I-E_+)x \in B_2^{1/2}\N(E_+) \subseteq B_2^{1/2}(B_2^{-1}(\St_+^{\perp}))=\Rt(B_2^{1/2}) \cap B_2^{-1/2}(\St_+^{\perp}) \subseteq \M_2^{\perp},$ where we used the fact that $E_+$ is $B_2$-symmetric and Proposition \ref{EBself}. Therefore, $P_{\M_2}B_2^{1/2} =B_2^{1/2}E_+$ in $\DD(E_+).$ Then, by multiplying both sides of the last equation by $B_2^{1/2},$ 
it follows that $$B_2E_+ \subset B_2^{1/2}P_{\M_2}B_2^{1/2}.$$


In a similar way, it can be proved that $P_{\M_3}B_3^{1/2} =B_3^{1/2}E_-$ in $\DD(E_-)$ and $B_3E_- \subset B_3^{1/2}P_{\M_3}B_3^{1/2}.$ 

Therefore, if $x \in \DD(E)=\DD(E_+)=\DD(E_-),$ then 
$$BEx=BE_+x+BE_-x=B_2E_+x-B_3E_-x=(B_2^{1/2}P_{\M_2}B_2^{1/2}-B_3^{1/2}P_{\M_3}B_3^{1/2})x.$$
Define $S:=B_2^{1/2}P_{\M_2}B_2^{1/2}-B_3^{1/2}P_{\M_3}B_3^{1/2},$ then $S \in L(\HH)^s$ and 
$BE \subset S=S^* \subset (BE)^*.$ Therefore $\ol{BE}=(BE)^*=S.$ 
{\hfill\qed\vspace{1ex}}

\bigskip
Let $B \in L(\HH)^s,$ $\St \subseteq \HH$ be a closed subspace of $\HH$ and $E$ be a $B$-symmetric closed projection onto $\St.$ The following example shows that if $B$ is not $\St$-weakly complementable then $BE$ may not admit a bounded selfadjoint extension to $\HH.$
We use the fact that a symmetric operator $T$ admits an extension to $\HH$ if and only if $T^* \in L(\HH)^s.$ Also, if $T$ is a densely defined operator such that $T \subset T^* \subset W \in L(\HH)$ then $W^*=\ol{T}=T^*=W.$

\begin{example} As in Example \ref{ex1}, let $\St \subseteq \HH$ be a closed subspace such that $\dim(\St)=\dim(\St^{\perp})=\infty.$ Take $a \in L(\St)^+$ such that $\Rt(a)$ is not closed and $\ol{\Rt(a)}=\St,$ and take $b: \St^{\perp} \ra \St$ such that $b$ is invertible. 
Consider $B=\begin{bmatrix}
		a & b \\ 
		b^* & c \\
\end{bmatrix}$ for any $c \in L(\St^{\perp})^s.$	

Let $x:=a^{\dagger}b$ and  $E=\begin{bmatrix}
		I & x\\ 
		0 & 0 \\
\end{bmatrix}$ where $\DD(E)=\St \oplus b^{-1}(\Rt(a)).$ Since $a^{\dagger}$ is closed and $b$ is bounded, $x$ is closed. Then $E$ is a closed projection onto $\St.$

Moreover, $$ax=aa^{\dagger}b=b \mbox { in } \DD(x)=b^{-1}(\Rt(a))$$ or equivalently, $ax \subset b.$ Applying Proposition \ref{propCuasi}, $E$ is $B$-symmetric.

By \cite[Corollary 2]{Holland}, $x^*=(a^{\dagger}b)^*=b^*a^{\dagger}$
and then $\DD(x^*)=\DD(a^{\dagger})=\Rt(a).$ Also, $E^*=\begin{bmatrix}
	I & 0\\ 
	x^* & 0 \\
\end{bmatrix}$  where $\DD(E^*)=\DD(x^*) \oplus \St^{\perp}=\Rt(a) \oplus \St^{\perp}.$ Let us see that $\DD(E^*B) \subsetneq \HH.$ In fact, let $h \in \St \setminus \Rt(a).$ Then $h=bk$ for some $k \in \St^{\perp},$ because $b$ is surjective.

Then $$Bk=\begin{bmatrix}
	a & b \\ 
	b^* & c \\
\end{bmatrix}  \begin{bmatrix}
0 \\ 
k \\
\end{bmatrix} =\begin{bmatrix}
h \\ 
ck \\
\end{bmatrix} \not \in \DD(E^*).$$

Therefore $k \not \in \DD(E^*B)$ and then $\DD(E^*B) \subsetneq \HH.$ Hence, by the comments above, neither $BE$ nor $E^*B$ admit a bounded selfadjoint extension to $\HH.$
\end{example}

\begin{prop} \label{PropBSym} Let $B \in L(\HH)^s$ and $\St \subseteq \HH$ be a closed subspace such that $B$ is $\St$-weakly complementable and the matrix decomposition of $B$ is as in \eqref{matrixW}. If $E$ is a closed  $B$-symmetric projection onto $\St$ then $E^* \in \mc{P}^*(B,\St).$
\end{prop}

\begin{dem} If $E$ is a closed $B$-symmetric projection onto $\St$ then $E^*$ is a densely defined closed projection with nullspace $\St^{\perp}$ and, by Proposition \ref{EBselfC2},  $E^*B=(BE)^* \in L(\HH)^s.$
	
On the other hand, by Proposition \ref{propaxb}, if the matrix decomposition of $E$ is as in \eqref{matrixE}, then $ax \subset b.$ Since $B$ is $\St$-weakly complementable, $b=\vert a \vert^{1/2}f,$ with $f$ the reduced solution of the equation $b=\vert a \vert^{1/2} h.$
Then, if $a=u\vert a\vert$ is the polar decomposition of $a,$  $(\vert a \vert^{1/2}f)z=(\vert a \vert^{1/2}  \vert a \vert^{1/2} u x) z \mbox{ for every } z \in \DD(x).$ Then $(\vert a \vert^{1/2} u x-f) \in \ol{\Rt(\vert a \vert^{1/2})} \cap \N(\vert a \vert^{1/2})=\{0\},$ so that $(\vert a \vert^{1/2} x) z= (uf) z \mbox{ for every } z \in \DD(x).$ Therefore $\vert a \vert^{1/2} x \subset uf.$ Then $(uf)^* \subset (\vert a \vert^{1/2} x)^*=x^* \vert a \vert^{1/2}$ and $x^* \vert a \vert^{1/2} \in L(\St, \St^{\perp}).$ Hence
$\Rt(\vert a \vert^{1/2}) \subseteq \DD(E^*).$
\end{dem}

\begin{prop} \label{PropBSym2} Let $B \in L(\HH)^s$  and $\St \subseteq \HH$ be a closed subspace such that $B\St \subseteq \St + (B\St)^{\perp}.$ Then $B$ is $\St$-weakly complementable and the pair $(B,\St)$ is quasi-complementable.
\end{prop}

\begin{dem} Since $\HH=\ol{B\St} \oplus (B\St)^{\perp} \subseteq \ol{\St + (B\St)^{\perp}},$  then $\St^{\perp} \cap \ol{B\St}=\{0\}.$ By Proposition \ref{EBselfC4}, the pair $(B,\St)$  is quasi-complementable, or equivalently, there exists a closed $B$-symmetric projection
$E$ onto $\St.$ If $B$ has matrix decomposition as in \eqref{matrixW} and $E=\begin{bmatrix}
	I & x \\ 
	0 & 0 \\
\end{bmatrix},$ where $x: \DD(x) \subseteq \St^{\perp} \ra \St$ is a densely defined closed linear operator then  $b^*=x^*a,$ by Proposition \ref{propaxb}.
 
To show that $B$ is $\St$-weakly complementable, we estimate $\vert P_{\St}-P_{{(B\St)}^{\perp}}\vert.$ From 
$$(P_{\St}-P_{{(B\St)}^{\perp}})^2=P_{\St}P_{\ol{B\St}}P_{\St}+P_{\St^{\perp}}P_{(B\St)^{\perp}}P_{\St^{\perp}},$$
it follows that $$ \vert P_{\St}-P_{{(B\St)}^{\perp}}\vert= \vert P_{\ol{B\St}} - P_{\St^{\perp}}\vert=(P_{\St}P_{\ol{B\St}}P_{\St})^{1/2}+ (P_{\St^{\perp}}P_{(B\St)^{\perp}}P_{\St^{\perp}})^{1/2}.
$$

Set $\Gamma_{E^*}:=(P_{\St^{\perp}}+P_{\N^{\perp}})^{1/2}.$ Then, 
$$\Rt((P_{\St}P_{\ol{B\St}}P_{\St})^{1/2}) \subseteq  \Rt(\vert P_{\ol{B\St}}-P_{\St^{\perp}}\vert)=\Rt( P_{\ol{B\St}}-P_{\St^{\perp}}) \subseteq \ol{B\St}+\St^{\perp} \subseteq \N^{\perp} + \St^{\perp} = \Rt(\Gamma_{E^*}).$$
On the other hand $\Rt(\vert P_{\St}BP_{\St}\vert)=\Rt(P_{\St}BP_{\St}) \subseteq \Rt(P_{\St}P_{\ol{B\St}} P_{\St}).$ In fact, if $y \in \Rt(P_{\St}BP_{\St})$ then there exists $x \in \HH$ such that $y=P_{\St}BP_{\St}x=P_{\St}P_{\ol{B\St}} (BP_{\St}x).$ Since $B\St \subseteq \St + (B\St)^{\perp},$ there exists $s \in \St$ and $t \in (B\St)^{\perp}$ such that $BP_{\St}x=s+t.$ Then $y=P_{\St}P_{\ol{B\St}} (s+t)=P_{\St}P_{\ol{B\St}}P_{\St} s.$ Therefore, by Douglas' Lemma and the monotonicity of the square root (see \cite{Pedersen}), there exists $c \geq 0$ such that
$$\vert P_{\St}BP_{\St}\vert \leq c^{1/2}  (P_{\St}P_{\ol{B\St}} P_{\St}),$$ or, equivalently,
$$\Rt(\vert a \vert^{1/2})=\Rt(\vert P_{\St}BP_{\St}\vert^{1/2}) \subseteq \Rt((P_{\St}P_{\ol{B\St}} P_{\St})^{1/2})\subseteq \Rt(\Gamma_{E^*})=\DD(E^*)=\DD(x^*) \oplus \St^{\perp}.$$ 
Then  $\Rt(\vert a \vert^{1/2}) \subseteq \DD(x^*).$ Applying Corollary \ref{corEBself}, $E^*B \in L(\HH)^s.$
Finally, since $E^*$ is a closed densely defined projection with $\N(E^*)=\St^{\perp},$ $E^*B \in L(\HH)^s$ and $\Rt(\vert a \vert^{1/2}) \subseteq \DD(E^*),$ by Theorem \ref{SWC}, $B$ is $\St$-weakly complementable.
\end{dem}

\begin{cor} Let $B \in L(\HH)^s$  and $\St \subseteq \HH$ be a closed subspace. Suppose that the matrix decomposition of $B$ is as in \eqref{matrixW}. If the pair $(B,\St)$ is quasi-complementable and $\Rt(\vert a \vert^{1/2}) \subseteq \Rt(P_{\St}P_{\ol{B\St}})$ then $B$ is $\St$-weakly complementable.
\end{cor}

\begin{dem} Observe that $\Rt(\vert a \vert^{1/2}) \subseteq \Rt(P_{\St}P_{\ol{B\St}})$ if and only if $\Rt(\vert a \vert^{1/2}) \subseteq \Rt((P_{\St}P_{\ol{B\St}} P_{\St})^{1/2}).$ Let $E=P_{\St \pl \N},$ with $\N\subseteq (B\St)^{\perp}.$ Then, if $\Rt(\vert a \vert^{1/2}) \subseteq \Rt((P_{\St}P_{\ol{B\St}} P_{\St})^{1/2}),$ by using arguments similar to those found in the proof of Proposition  \ref{PropBSym2} we get that $E^*$ is a closed densely defined projection with $\N(E^*)=\St^{\perp},$ $E^*B \in L(\HH)^s$ and $\Rt(\vert a \vert^{1/2}) \subseteq \DD(E^*).$ Then, by Theorem \ref{SWC}, $B$ is $\St$-weakly complementable.
\end{dem}

\subsection{Applications: Schur complements of selfadjoint operators}
We recall the definition of Schur complement for an $\St$-weakly complementable selfadjoint operator.

\begin{Def}  Let $B \in L(\HH)^s$ and $\St \subseteq \HH$ be a closed subspace such that $B$ is $\St$-weakly complementable. When $B$ is as in \eqref{matrixW}, let $f$ be the reduced solution of $b=\vert a\vert^{1/2}x$ and  $a=u \vert a \vert$ the polar decomposition of $a.$ The \emph{Schur complement} of $B$ to $\St$ is defined as
	$$B_{ / \St}:= \begin{bmatrix}
	0 & 0 \\ 
	0 & c-f^*uf \\
	\end{bmatrix}$$ 
and	$B_{\St} := B- B_{/ \St}$ is the \emph{compression} of $B$ to $\St.$
\end{Def}

When $B \in L(\HH)^+$ this formula gives the usual Schur complement, see \cite[Theorem 3]{Shorted2}.

\bigskip
When the operator $B$ is $\St$-complementable, the Schur complement can be written as $B_{/ \St}=(I-F)B,$ for any bounded projection with $\N(F)=\St^{\perp}$ such that $(FB)^*=FB.$ In fact, from \cite[Corollary 3.12]{Contino4} it suffices to take $F=Q^*,$ for any  $Q \in \mc{P}(B,\St).$

A similar formula for $B_{/ \St}$ can be given when $B$ is $\St$-weakly complementable. In this case the projection need not be bounded, but it is a semiclosed densely defined projection with closed nullspace.

\begin{thm}[cf. {\cite[Theorem 3.14]{Contino4}}] \label{thmwc} Let $B \in L(\HH)^s$ and $\St$ be a closed subspace of $\HH$. Suppose that $B$ is $\St$-weakly complementable then
$$B_{ / \St} = (I-E)B,$$ for every $E \in \mc{P}^*(B,\St).$ 
\end{thm}

\begin{dem} Let $E \in \mc{P}^*(B,\St).$ Suppose that the matrix decomposition of $E$ is as in \eqref{ESCDecomp} and that of $B$ is as in \eqref{matrixW}. Then, by Lemma \ref{ProjectionBself}, $ya=b^*$ and since $\Rt(\vert a \vert^{1/2}) \subseteq \DD(E),$ by  Corollary \ref{corKaufman2}, $E \vert a\vert^{1/2} \in L(\St, \HH),$ so that $y \vert a\vert^{1/2} \in L(\St, \St^{\perp})$. Hence 
	$$(I-E)B=\begin{bmatrix} 
	0 & 0 \\ 
	-y  & I\\  
	\end{bmatrix}\begin{bmatrix} 
	a & b \\ 
	b^*  & c\\  
	\end{bmatrix}=\begin{bmatrix} 
	0 & 0 \\ 
	-ya+b^*  & c-yb	\end{bmatrix}=\begin{bmatrix} 
	0 & 0 \\ 
	0  & c-yb\\  
	\end{bmatrix}.$$ Let $f$ be the reduced solution of $b=\vert a\vert^{1/2}x$ and  $a=u \vert a \vert$ the polar decomposition of $a.$ Then $yb=f^*uf.$ In fact, since $ya=b^*$ we have that 
	$y\vert a \vert=f^* \vert a \vert^{1/2}u=f^*u \vert a\vert^{1/2}.$ Then $y\vert a \vert^{1/2}=f^*u$ on $\Rt(\vert a \vert^{1/2}),$ and since $y \vert a \vert^{1/2}$ is bounded, $y\vert a \vert^{1/2}=f^*u$ on $\ol{\Rt(\vert a \vert^{1/2})}.$ Then $yb=y\vert a \vert^{1/2}f=f^*uf$ because $\Rt(f)\subseteq \ol{\Rt(\vert a \vert^{1/2})}.$ Hence $(I-E)B=B_{ / \St}.$
\end{dem}

In particular, Theorem \ref{thmwc} gives a formula for the Schur complement of any positive operator $B$ to $\St$ in terms of semiclosed projections. 

\bigskip
Different definitions where given for the minus order, for example, using generalized inverses in the matrix case, see \cite{Mit86}. We give the following definition, equivalent to those appearing in \cite{Sem10} and \cite{Minus}.

\begin{Def} Let $A, B\in L(\HH)$, we write $A \minus B$  if there exist projections $P, Q \in \Q$  such that $A=PB$ and $A^*=QB^*$.
\end{Def}

It was proved in \cite{Sem10} and \cite{Minus} that $\minus$ is a partial order, known as the minus order for operators. In \cite{AntCorSto06}, it was shown that  $A\minus B$ if and only if the sets $ \ol{\Rt(A)}\dot{+} \ol{\Rt(B-A)} \mbox{ and } \ol{\Rt(A^*)}\dot{+}\ol{\Rt(B^*-A^*)}$ are closed.
In \cite[Theorem 3.3]{Minus},  another characterization of the minus order in terms of the range additivity property was given:
$$A\minus B \mbox{ if and only if } \Rt(B)=\Rt(A)\dot{+}\Rt(B-A) \mbox{ and } \Rt(B^*)=\Rt(A^*)\dot{+}\Rt(B^*-A^*).$$

In view of this equivalence, the left minus order  was defined in \cite{Minus} for operators in $L(\HH).$ This notion is weaker than the minus order, in the infinite dimensional setting. 

\begin{Def} Let $A, B\in L(\HH)$, we write  $A \lminus B$  if $\Rt(B)=\Rt(A)\dot{+}\Rt(B-A)$.
\end{Def}

The relation $\lminus$ is a partial order, see \cite{Minus}. The following is a characterization of the left minus order in terms of a semiclosed projection.
\begin{prop}[cf. {\cite[Proposition 3.13]{Minus}}]\label{leftminus and projections} Let $A, B\in L(\HH)$. Then $A \lminus B$ if and only if there exists $P \in \Q_{SC}$ such that $A=PB$ and $\Rt(A)\subseteq \Rt(B)$.
\end{prop}

\begin{dem} Suppose that $A \lminus B$ then $\Rt(B)=\Rt(A)\dot{+}\Rt(B-A)$ so that $\Rt(A)\subseteq \Rt(B)$. Define $P=P_{\Rt(A) \pl {\Rt(B-A)\oplus \N(B^*)}}$. Then $P \in \Q_{SC}$ and it is easy to check that $A=PB$. 
	
Conversely, suppose that $A=PB$ for $P \in \Q_{SC}$ and $\Rt(A)\subseteq \Rt(B).$ Then $\Rt(B)=\Rt(A)+\Rt(B-A)$ and the sum is direct because $\Rt(A)\subseteq \Rt(P)$ and $\Rt(B-A)\subseteq \N(P).$
\end{dem}

A different generalization of the minus order was introduced by Arias et al. in \cite{Arias2}:

\begin{Def} Given $A,B \in L(\HH),$ we  write $A \prec B$ if there exist  $Q, P \in \Q_{SC}$ with closed ranges such that $A=QB$ and $A^*=PB^*.$
\end{Def}

The relation $\prec$ is a partial order in $L(\HH)$ (see \cite[Lemma 4.5]{Arias2}).
In \cite[Lemma 4.4]{Arias2}, it was proved that $A \prec B$ if and only if $\ol{\Rt(A)} \cap \Rt(B-A)=\{0\}$ and $\ol{\Rt(A^*)} \cap \Rt(B^*-A^*)=\{0\}.$ 

\bigskip
Let $B \in L(\HH)^s$ and $\St$ be a closed subspace of $\HH$. Suppose that the matrix decomposition of $B$ is as in \eqref{matrixW}. Denote by $\tilde{\Q}:=\{Q\in  \Q_{SC}:  \Rt(Q) \mbox{ is closed}  \mbox{ and } \Rt(\vert a \vert^{1/2}) \oplus \St^{\perp} \subseteq \DD(Q)\}.$ Note that if $E \in \mc{P}^*(B,\St)$ then $I-E \in \tilde{Q}.$ 
Define
$$\mc{M}(B,\St):=\{X \in L(\HH)^s : \Rt(X) \subseteq \St^{\perp} \mbox{ and } \ X=QB,  \mbox{ for some } Q \in \tilde{Q}\}.$$

\begin{thm} Let $B \in L(\HH)^s$ and $\St$ be a closed subspace of $\HH$ such that $B$ is $\St$-weakly complementable. Then
	$$B_{ / \St}= \max_{\prec} \ \mc{M}(B,\St).$$
\end{thm}

\begin{dem}  Let $E \in \mc{P}^*(B,\St).$ Then, by Theorem \ref{thmwc}, $B_{ / \St} = (I-E)B \in L(\HH)^s$ and $\Rt(B_{ / \St}) \subseteq \St^{\perp}.$ Set $Q:=I-E,$ then $Q \in \tilde{Q}$ and $B_{ / \St} = QB.$ So that $B_{ / \St}  \in \mc{M}(B,\St).$ 

On the other hand,  let $X \in \mc{M}(B,\St).$ Then $X \in L(\HH)^s,$ $\Rt(X) \subseteq \St^{\perp}$ and there exists $Q \in \tilde{Q}$ such that $X=QB. $ Suppose that the matrix decomposition of $B$ is as in \eqref{matrixW} and let $f$ be the reduced solution of $b=\vert a\vert^{1/2}x$ and  $a=u \vert a \vert$ the polar decomposition of $a.$ Then 

$$B=\begin{bmatrix} 
a & \vert a\vert^{1/2}f \\ 
f^*\vert a\vert^{1/2}  & c\\  
\end{bmatrix}=\Gamma U \Gamma,$$

where $\Gamma:=\begin{bmatrix} 
\vert a \vert^{1/2} & 0 \\ 
0 & I\\  
\end{bmatrix} \in L(\HH)^+$ and $U:=\begin{bmatrix} 
u & f \\ 
f^* & c\\  
\end{bmatrix} \in L(\HH)^s.$  

Since $Q$ is semiclosed  and $\Rt(\Gamma)=\Rt(\vert a \vert^{1/2}) \oplus \St^{\perp} \subseteq \DD(Q)$, by Corollary \ref{corKaufman2}, $Q \Gamma \in L(\HH).$  Therefore
$$X=QB=Q\Gamma U \Gamma =X^*=(Q\Gamma U \Gamma )^*= \Gamma U (Q\Gamma)^*.$$ 
Let $E \in \mc{P}^*(B,\St).$ Since $\Rt(\Gamma)\subseteq \DD(E)=\DD(I-E)$, by Corollary \ref{corKaufman2}, $(I-E) \Gamma \in L(\HH).$ Also,
\begin{align*}
X&=(I-E)X=(I-E)\Gamma U (Q\Gamma)^*=X^*=((I-E)\Gamma U (Q\Gamma)^*)^*\\
&=Q\Gamma U ((I-E)\Gamma)^*=Q((I-E)\Gamma U \Gamma)^*=Q((I-E)B)^*=QB_{ / \St},
\end{align*}
where we used Theorem \ref{thmwc}. Then $X \prec B_{ / \St}.$ 
\end{dem}

\section*{Acknowledgements}
We thank the anonymous referee for carefully reading our manuscript and helping us to improve it with several useful comments.

Maximiliano Contino was supported by the UBA's Strategic Research Fund 2018 and CONICET PIP 0168. Alejandra Maestripieri was supported by CONICET PIP 0168.  The work of Stefania Marcantognini was done during her stay at  the Instituto Argentino de Matem\'atica with an appointment funded by the CONICET. She is greatly grateful to the institute for its hospitality and to the CONICET for financing her post.


\end{document}